\documentclass[11pt]{article}
\usepackage[a4paper,margin=1in]{geometry}
\usepackage{amsmath,amssymb,amsfonts,amsthm,mathtools}
\usepackage{lmodern}
\usepackage{microtype}
\usepackage{graphicx}
\graphicspath{{numerics/}{../../output/section5_en_a20_eps9_20261009/figures/}{../../empirical_section6_histdata/english_20261011/figures/}}
\usepackage{booktabs}
\usepackage{float}
\usepackage{xcolor}
\usepackage{hyperref}

\definecolor{primary}{RGB}{24,64,122}
\hypersetup{
  colorlinks=true,
  linkcolor=primary,
  citecolor=primary,
  urlcolor=primary
}

\newcommand{\R}{\mathbb{R}}
\newcommand{\E}{\mathbb{E}}
\newcommand{\Prob}{\mathbb{P}}
\newcommand{\dd}{\,\mathrm{d}}
\newcommand{\cF}{\mathcal{F}}
\newcommand{\cL}{\mathcal{L}}
\newcommand{\cB}{\mathcal{B}}
\newcommand{\cP}{\mathcal{P}}

\newcommand{\Lip}{\operatorname{Lip}}
\newcommand{\TV}{\mathrm{TV}}
\newcommand{\BL}{\mathrm{BL}}
\newcommand{\HS}{\mathrm{HS}}
\newcommand{\op}{\mathrm{op}}
\newcommand{\ProjX}{\operatorname{Proj}_X}

\newtheorem{assumption}{Assumption}
\newtheorem{theorem}{Theorem}[section]
\newtheorem{proposition}[theorem]{Proposition}
\newtheorem{lemma}[theorem]{Lemma}

\theoremstyle{definition}

\theoremstyle{remark}
\newtheorem{remark}[theorem]{Remark}

\numberwithin{equation}{section}

\title{Dimension Reduction in Nonlinear Filtering of Multiscale\\
Stochastic Volatility from Block Realized Measures}
\author{Qi Zhang\thanks{Corresponding author. Beijing Institute of Mathematical Sciences and Applications, Beijing 101408, China. Email: \texttt{qzhang@bimsa.cn}}, Yiou Zhou\thanks{ Beijing Institute of Mathematical Sciences and Applications, Beijing 101408, China. Email: \texttt{ minorca1028@outlook.com}}}

\date{\today}

\begin{document}
\maketitle

\begin{abstract}
We study nonlinear filtering for multivariate slow--fast stochastic volatility diffusions observed through blockwise logarithmic realized variance proxies under a Gaussian observation model. Because each observation depends nonlinearly on an entire trajectory segment, exact filtering requires a high-dimensional path-integral recursion that is stiff at the fast time scale. We prove that the slow marginal of the full filter converges on finite horizons, in expected bounded-Lipschitz distance, to an autonomous reduced filter as the block length vanishes while remaining large relative to the fast scale. We also derive quantitative localized error bounds revealing a tradeoff between fast ergodic averaging and slow-state freezing, and identify the state-dependent convexity adjustment induced by averaging variance prior to taking logarithms. We further establish convergence of the reduced filter to its continuous-time weak Zakai limit and confirm these mechanisms numerically.
\end{abstract}



\section{Introduction}
\label{sec:introduction}

Asset return volatility exhibits pronounced multiscale dynamics, combining persistent macroeconomic or regime factors with rapidly mean-reverting intraday fluctuations \cite{AndersenBollerslev1998,BarndorffNielsenShephard2002,FouquePapanicolaouSircar2000,FouquePapanicolaouSircarSolna2011}. Because instantaneous spot variance is latent, high-frequency financial econometrics constructs nonparametric proxies by aggregating intraday squared returns over discrete time windows, yielding realized variance estimators that approximate integrated variance over each block \cite{AndersenBollerslev1998,JacodProtter2012}. Filtering the persistent slow volatility factors from such blockwise logarithmic realized measures poses a nonstandard multiscale inference problem: each observation depends nonlinearly on an entire continuous trajectory segment rather than a single terminal state, microscopic fast-scale stiffness makes direct high-dimensional filtering computationally prohibitive, and averaging spot variance prior to taking logarithms induces a systematic convexity shift that must be captured in any reduced model.

Consider a multivariate stochastic volatility system for $d$ assets driven by $m$ slow volatility factors $X_t^\varepsilon \in \R^m$ and $n$ fast volatility factors $Y_t^\varepsilon \in \R^n$. For fixed loading vectors $\ell_r \in \R^m$ and $c_r \in \R^n$, the instantaneous spot variance of asset $r \in \{1,\dots,d\}$ is
\[
  V_{t,r}^\varepsilon = v_r(X_t^\varepsilon, Y_t^\varepsilon),
  \qquad
  v_r(x,y) := \exp(\ell_r^\top x + c_r^\top y).
\]
Measured in units of the slow characteristic time $\tau_s$, the latent state $(X_t^\varepsilon, Y_t^\varepsilon) \in \R^m \times \R^n$ evolves according to the slow--fast stochastic differential equation (SDE)
\begin{align}
  \dd X_t^\varepsilon
  &= f(X_t^\varepsilon, Y_t^\varepsilon)\dd t
   + \sigma_X(X_t^\varepsilon, Y_t^\varepsilon)\dd W_t^X,
  \label{eq:slow}\\
  \dd Y_t^\varepsilon
  &= \varepsilon^{-1} b(X_t^\varepsilon, Y_t^\varepsilon)\dd t
   + \varepsilon^{-1/2} \tau(X_t^\varepsilon, Y_t^\varepsilon)\dd W_t^Y,
  \label{eq:fast}
\end{align}
where $\varepsilon := \tau_f / \tau_s \ll 1$ is the dimensionless scale-separation parameter, and $W^X, W^Y$ are independent standard Brownian motions in $\R^m$ and $\R^n$, respectively.

Partitioning time into observation blocks $I_k := ((k-1)\delta, k\delta]$ of normalized length $\delta = \Delta/\tau_s$ (with physical window length $\Delta$ and grid points $t_k := k\delta$), define the block-averaged spot variance and its componentwise logarithm for $r = 1,\dots,d$ by
\begin{equation}
  A_{k,r}^\varepsilon := \frac{1}{\delta}\int_{I_k} v_r(X_t^\varepsilon, Y_t^\varepsilon)\dd t,
  \qquad
  H_{k,r}^\varepsilon := \log A_{k,r}^\varepsilon,
  \label{eq:block-functional}
\end{equation}
and collect these entries in $A_k^\varepsilon, H_k^\varepsilon \in \R^d$. As reviewed in Appendix~\ref{app:rv-clt}, a fixed-$\varepsilon$ high-frequency central limit expansion for realized variance motivates an additive Gaussian error model for logarithmic realized variance proxies, while full multivariate high-frequency asymptotics with path-dependent covariances and microstructure noise require separate econometric machinery \cite{BarndorffNielsenShephard2004,BarndorffNielsenEtAl2008}. We therefore study the vector Gaussian block observation model
\begin{equation}
  Z_k = H_k^\varepsilon + R_\delta^{1/2}\xi_k,
  \qquad \xi_k \overset{\mathrm{i.i.d.}}{\sim} \mathcal{N}(0, I_d),
  \qquad Z_k \in \R^d,
  \label{eq:realized-measure-observation}
\end{equation}
where $(\xi_k)_{k\ge 1}$ are independent of $(X^\varepsilon, Y^\varepsilon)$ and the deterministic covariance matrix $R_\delta = R_0/\delta + \widetilde{R}_\delta$ satisfies $R_0 \succ 0$ with a bounded positive semidefinite remainder $\widetilde{R}_\delta$ (specified in \eqref{eq:observation-noise-scaling}). Throughout, \eqref{eq:realized-measure-observation} serves as an idealized Gaussian observation model that isolates the filtering effects of multiscale block averaging from high-frequency microstructure and quarticity effects.

The coupled system \eqref{eq:slow}--\eqref{eq:fast} and \eqref{eq:realized-measure-observation} departs from classical continuous-time or discrete point-observation nonlinear filtering \cite{BainCrisan2009} in three fundamental respects:
\begin{enumerate}
  \item[\textup{(1)}] \emph{Path-dependent nonlinear observables.} Each observation $Z_k$ depends on the continuous trajectory segment $\mathcal{X}_k^\varepsilon := (X_s^\varepsilon, Y_s^\varepsilon)_{s\in I_k}$ through the nonlinear integral functional $H_k^\varepsilon$.
  \item[\textup{(2)}] \emph{Path-integral Bayesian updates.} Conditional on $\mathcal{X}_k^\varepsilon$, the observation $Z_k$ has a Gaussian density centered at $H_k^\varepsilon$, so Bayes updating requires a joint path-integral transition operator rather than pointwise multiplication of a terminal-state likelihood.
  \item[\textup{(3)}] \emph{Double asymptotic scale window.} The block length $\delta$ lies between the fast and slow characteristic scales:
  \begin{equation}
    \varepsilon \ll \delta \ll 1,
    \qquad\text{equivalently}\qquad
    \tau_f \ll \Delta \ll \tau_s.
    \label{eq:scale-window}
  \end{equation}
  Within each block $I_k$, the fast factor completes $O(\delta/\varepsilon) \gg 1$ ergodic mixing cycles, while the slow factor moves by only $O(\sqrt{\delta}) \ll 1$.
\end{enumerate}

The exact Bayesian benchmark is the $(m+n)$-dimensional conditional law $\Pi_k^\varepsilon := \Prob((X_{t_k}^\varepsilon, Y_{t_k}^\varepsilon) \in \cdot \mid \cF_k^Z)$, with its slow marginal $\pi_k^\varepsilon := \Pi_k^\varepsilon(\cdot \times \R^n)$ as our primary inference target. When $\varepsilon \ll 1$, direct numerical implementation of $\Pi_k^\varepsilon$ is computationally prohibitive because the fast drift and diffusion scale as $\varepsilon^{-1}$ and $\varepsilon^{-1/2}$, forcing explicit simulation or particle filters to resolve the microscopic scale $\varepsilon$. Moreover, rapid ergodicity of the frozen fast diffusion does not by itself yield a naive factorization of the posterior, since the moving slow state and the path-integral likelihood couple slow and fast trajectories across every block.

In the scale window \eqref{eq:scale-window}, averaging the slow coefficients against the frozen fast invariant measure $\mu_x$ yields an autonomous $m$-dimensional homogenized signal $\bar X_t$ with drift $\bar f(x) = \int_{\R^n} f(x,y)\,\mu_x(\dd y)$ and diffusion covariance $\bar a(x) = \int_{\R^n} \sigma_X(x,y)\sigma_X(x,y)^\top\,\mu_x(\dd y)$. Because each block averages the exponential spot variance \emph{before} taking the logarithm, the effective observation map has components
\begin{equation}
  \bar H_r(x)
  := \log \int_{\R^n} v_r(x,y)\,\mu_x(\dd y)
  = \ell_r^\top x + \log \int_{\R^n} e^{c_r^\top y}\,\mu_x(\dd y),
  \qquad 1 \le r \le d.
  \label{eq:homogenized-realized-measure-map}
\end{equation}
By Jensen's inequality, whenever $\int_{\R^n} c_r^\top y\,\mu_x(\dd y) = 0$ and $c_r^\top Y$ is nonconstant under $\mu_x$, the logarithmic moment-generating term in \eqref{eq:homogenized-realized-measure-map} introduces a strictly positive \emph{convexity adjustment}. For example, in the scalar centered Ornstein--Uhlenbeck (OU) model with $\mu_x = \mathcal{N}(0,\nu_f^2/(2\kappa_f))$ and unit loadings,
\begin{equation}
  \bar H(x) = x + \frac{\nu_f^2}{4\kappa_f}.
  \label{eq:ou-jensen-correction}
\end{equation}
More generally, when $\mu_x$ depends on $x$, this correction modifies both the level and gradient of $\bar H$, and omitting it induces a systematic upward bias in the inferred slow volatility state.

This paper quantifies how accurately the autonomous $m$-dimensional reduced filter $\bar\pi_k$ (built from $\bar X_t$ and the terminal-state likelihood centered at $\bar H$) approximates the slow marginal $\pi_k^\varepsilon$ uniformly over $k \le T/\delta$ blocks, and characterizes how the error scales jointly with $\varepsilon$ and $\delta$. On a single observation block, replacing the path functional $H_k^\varepsilon$ by the endpoint map $\bar H(X_{t_k}^\varepsilon)$ incurs two competing physical errors governed by
\begin{equation}
  r_{\varepsilon,\delta} := \sqrt{\frac{\varepsilon}{\delta}} + \sqrt{\delta},
  \label{eq:leading-error-intro}
\end{equation}
where $\sqrt{\varepsilon/\delta}$ measures incomplete fast ergodic averaging over $I_k$ and $\sqrt{\delta}$ measures freezing the slow trajectory at the block endpoint. Balancing these two contributions selects the geometric-mean window scaling $\delta \asymp \sqrt{\varepsilon}$, or $\Delta^* \asymp \sqrt{\tau_f \tau_s}$ in physical time.

Under primitive dissipativity, ellipticity, regularity, and moment conditions (Assumption~\ref{ass:dissipative}), our main convergence result (Theorem~\ref{thm:filter-rate}) proves that for every finite horizon $T < \infty$, every sufficiently large localization radius $R$, and an integer Sobolev index $s > 4 + m/2$, the expected bounded-Lipschitz posterior discrepancy satisfies
\[
  \max_{k\delta \le T} \E^{\Prob}\, d_{\BL}(\pi_k^\varepsilon, \bar\pi_k)
  \le C_{T,R}\!\left[r_{\varepsilon,\delta} + \left(\frac{\varepsilon}{\delta}\right)^{1/s}\right] + O(R^{-q})
  \le C_{T,R}\, r_{\varepsilon,\delta}^{1/s} + O(R^{-q}),
\]
where the explicit polynomial tail $O(R^{-q})$ bounds the probability that the slow state exits the ball of radius $R$ before time $T$. Letting $\varepsilon/\delta \to 0$ and $\delta \to 0$ at fixed $R$, followed by $R \to \infty$, establishes finite-horizon convergence of the reduced filter to the slow marginal of the full filter. In addition, Theorem~\ref{thm:discrete-continuous-filter} proves that as $\delta \downarrow 0$, the discrete reduced filter driven by its own homogenized observations converges in expected total variation at localized rate $O(\sqrt{\delta})$ to the continuous-time nonlinear filter governed by the weak Zakai stochastic partial differential equation (SPDE) \cite{YY1998,Zakai1969,LuoYau2013}.

To avoid the random normalizing constants in Bayes' rule, our proof of Theorem~\ref{thm:filter-rate} operates on unnormalized filtering measures under a common product Gaussian reference measure \cite{KallianpurStriebel1968,BainCrisan2009} and advances through three key steps:
\begin{enumerate}
  \item[\textup{(i)}] \emph{Localization and parameter-regular Poisson correctors.} Applying a smooth radial truncation to the slow coefficient argument, we derive Gaussian-weighted $C_b^{s+3}$ parameter regularity for the frozen Poisson correctors and invariant averages directly from the primitive coefficients, without assuming differentiability of the invariant density (Lemma~\ref{lem:cutoff-poisson}).
  \item[\textup{(ii)}] \emph{Multiblock observation reduction.} Combining Poisson-corrector It\^o expansions with Gaussian Hellinger comparison, we replace the path-dependent block observations by endpoint observations in total variation at cost $O(r_{\varepsilon,\delta})$; the $O(\delta)$ Fisher information per block precisely cancels the $O(1/\delta)$ accumulation across blocks (Lemma~\ref{lem:block-endpoint-experiment}).
  \item[\textup{(iii)}] \emph{Weighted Sobolev transition telescoping.} Telescoping the full and averaged transition kernels in polynomially weighted Sobolev spaces yields an $O(\varepsilon/\delta)$ dual Sobolev error and an $O((\varepsilon/\delta)^{1/s})$ bounded-Lipschitz error (Proposition~\ref{prop:forward-adapted-energy}). A Bayes quotient transfer lemma and pathwise coupling then remove the normalization and localization, respectively.
\end{enumerate}

Dimension reduction for multiscale nonlinear filtering has been investigated in continuous time by Park, Sowers, and Namachchivaya \cite{ParkSowersNamachchivaya2010}, via homogenization by Imkeller, Namachchivaya, Perkowski, and Yeong \cite{ImkellerEtAl2013}, and for correlated slow--fast systems by Beeson, Namachchivaya, and Perkowski \cite{BeesonNamachchivayaPerkowski2022}. In discrete time, Papavasiliou \cite{Papavasiliou2007} studies particle filters built from averaged dynamics at fixed slow-time intervals, while Papanicolaou \cite{Papanicolaou2012} analyzes a regime-switching volatility model with a fast Ornstein--Uhlenbeck factor and path-dependent interval likelihoods at a fixed sampling interval and finite observation sequence. In contrast to these continuous-time or finite-state fixed-window settings, our analysis treats an $\R^m$-valued slow diffusion under nonlinear logarithmic block-average observations, explicitly identifies the state-dependent convexity-adjusted map \eqref{eq:homogenized-realized-measure-map}, and establishes quantitative posterior error bounds in the joint regime $\varepsilon/\delta \to 0$, $\delta \to 0$ across $O(1/\delta)$ updates.

On the stochastic analysis side, our signal averaging builds on Khasminskii \cite{Khasminskii1968} and Pavliotis and Stuart \cite{PavliotisStuart2008}, together with the noncompact Poisson-equation methodology of Pardoux and Veretennikov \cite{PardouxVeretennikov2001,PardouxVeretennikov2003,PardouxVeretennikov2005}. To derive parameter regularity of Poisson correctors and invariant averages from primitive coefficients, we adapt the transfer-identity framework of R\"ockner and Xie \cite{RocknerXie2021} from polynomially growing forcings to Gaussian-weighted H\"older spaces capable of handling exponential variance observables.

In financial econometrics, realized variance and multivariate realized covariation were developed by Andersen and Bollerslev \cite{AndersenBollerslev1998}, Barndorff-Nielsen and Shephard \cite{BarndorffNielsenShephard2002,BarndorffNielsenShephardCovariation2004}, and Jacod and Protter \cite{JacodProtter2012}, inspiring state-space volatility models with realized measures such as Realized GARCH \cite{HansenHuangShek2012} and realized stochastic volatility \cite{TakahashiOmoriWatanabe2009,KoopmanScharth2013,YamauchiOmori2020}. While these econometric frameworks motivate logarithmic realized measures as volatility observations, their high-frequency asymptotic laws involve path-dependent quarticity covariances and must be distinguished from the deterministic-covariance Gaussian block model analyzed here.

\paragraph{Notation and Conventions.}
Throughout the paper, the dimensions $m,n,d$ are fixed and $0 < \varepsilon \le \delta \le 1$, with $r_{\varepsilon,\delta} := \sqrt{\varepsilon/\delta} + \sqrt{\delta}$. For vectors, $|\cdot|$ denotes the Euclidean norm; for matrices, $\|\cdot\|_{\HS}$ denotes the Hilbert--Schmidt norm, $\|\cdot\|_{\op}$ the spectral (operator) norm, and $A:B := \operatorname{tr}(A^\top B)$ the Frobenius inner product. We write $a(x,y) := \sigma_X(x,y)\sigma_X(x,y)^\top$ and $A_F(x,y) := \tau(x,y)\tau(x,y)^\top$. Multi-index partial derivatives are denoted by $D_x^\alpha D_y^\gamma$. For $j \ge 0$, $C_b^j$ denotes functions with bounded continuous derivatives up to order $j$, and $C_{\mathrm{loc}}^j$ denotes $j$-times continuously differentiable functions.

For a bounded Lipschitz function $\varphi:\R^m \to \R$, its bounded-Lipschitz norm is
\[
  \|\varphi\|_{\BL} := \|\varphi\|_\infty + \Lip(\varphi),
  \qquad
  \Lip(\varphi) := \sup_{x \ne y}\frac{|\varphi(x)-\varphi(y)|}{|x-y|}.
\]
For finite positive Borel measures $\mu,\nu$ on $\R^m$, the bounded-Lipschitz distance \cite[Section~11.3]{Dudley2002} and the total variation norm of a finite signed measure $\zeta$ are defined by
\[
  d_{\BL}(\mu,\nu) := \sup_{\|\varphi\|_{\BL} \le 1} |\mu(\varphi) - \nu(\varphi)|,
  \qquad
  \|\zeta\|_{\TV} := \sup_{\|\varphi\|_\infty \le 1} |\zeta(\varphi)|.
\]
The metric $d_{\BL}$ metrizes weak convergence on $\cP(\R^m)$ and satisfies $|\mu(1)-\nu(1)| \le d_{\BL}(\mu,\nu) \le \|\mu-\nu\|_{\TV}$ for all finite positive measures $\mu,\nu$. We define the fast-variable Gaussian Lyapunov weight $V_\eta(y) := e^{\eta|y|^2}$ with $\eta > 0$. Generic positive constants $C$ may vary from line to line and may depend on $m,n,d$, the observation loadings, and primitive model constants; dependencies on $T$, the localization radius $R$, and the initial moment bound $M_0$ are indicated by subscripts.

\paragraph{Organization of the Paper.}
Section~\ref{sec:preliminaries} collects global notation, formulates the path-integral filter, states Assumption~\ref{ass:dissipative}, proves frozen ergodicity, and derives the convexity-adjusted reduced filter. Section~\ref{sec:convergence} develops the convergence theory: Subsection~\ref{subsec:main-theorem-and-scale-balancing} states Theorem~\ref{thm:filter-rate} and the unnormalized error decomposition; Subsection~\ref{subsec:localization-poisson} establishes localization and Poisson corrector regularity; Subsection~\ref{subsec:observation-reduction} performs the block-to-endpoint observation reduction; Subsection~\ref{subsec:transition-telescoping} executes the weighted Sobolev transition telescoping and completes the proof of Theorem~\ref{thm:filter-rate}; and Subsection~\ref{subsec:continuous-time-limit} proves the continuous-time weak Zakai limit (Theorem~\ref{thm:discrete-continuous-filter}). Section~\ref{sec:examples} presents controlled Monte Carlo experiments verifying balanced-scale convergence, computational speedup, the U-shaped window tradeoff, the convexity adjustment, and slow-state recovery. Section~\ref{sec:conclusion} concludes, and Appendix~\ref{app:rv-clt} reviews the high-frequency price-process expansion motivating \eqref{eq:realized-measure-observation}.

\section{Mathematical Setup and Multiscale Filtering Formulation}
\label{sec:preliminaries}

\subsection{Gaussian Block Observation Model and Full Path-Integral Filter}
\label{subsec:filtrations-and-experiment}

Recall the block-averaged spot variance $A_k^\varepsilon$ and log-variance observable $H_k^\varepsilon$ defined in \eqref{eq:block-functional} over $I_k = (t_{k-1}, t_k]$. Under the observation model \eqref{eq:realized-measure-observation}, the deterministic noise covariance satisfies
\begin{equation}
  R_\delta = \frac{R_0}{\delta} + \widetilde{R}_\delta,
  \qquad R_0 \succ 0,
  \qquad 0 \preceq \widetilde{R}_\delta \preceq \overline{R} I_d,
  \qquad r_0 := \lambda_{\min}(R_0) > 0.
  \label{eq:observation-noise-scaling}
\end{equation}
In particular, $\|R_\delta^{-1}\|_{\op} \le \delta/r_0$, so the Gaussian Fisher information for the observation mean is of order $O(\delta)$ per block and its sum over $k \le T/\delta$ blocks remains bounded on every finite horizon $T$. Off-diagonal entries of $R_\delta$ permit contemporaneous cross-asset observation noise; cross-asset realized covariances, which need not be positive, are excluded from the componentwise logarithmic map.

Let $\phi_\delta$ denote the centered $\mathcal{N}(0,R_\delta)$ probability density on $\R^d$:
\begin{equation}
  \phi_\delta(u) := (2\pi)^{-d/2}(\det R_\delta)^{-1/2}
  \exp\!\left(-\frac{1}{2} u^\top R_\delta^{-1} u\right),
  \qquad u \in \R^d.
  \label{eq:block-gaussian-density}
\end{equation}
Conditional on the trajectory segment $\mathcal{X}_k^\varepsilon := (X_t^\varepsilon, Y_t^\varepsilon)_{t\in I_k}$, the observation $Z_k$ is independent of past observations and has conditional density $\phi_\delta(z - H_k^\varepsilon)$.

Following the Kallianpur--Striebel reference-measure framework \cite{KallianpurStriebel1968,BainCrisan2009}, we introduce a reference probability measure $\mathbb{Q}_\delta$ under which the signal $(X^\varepsilon, Y^\varepsilon)$ retains its original law while the observations $(Z_j)_{j\ge 1}$ are i.i.d.\ $\mathcal{N}(0,R_\delta)$ vectors independent of the signal. Setting $\cF_T^{\mathrm{sig}} := \sigma((X_t^\varepsilon,Y_t^\varepsilon): 0 \le t \le T)$ and $\cF_k^Z := \sigma(Z_1,\dots,Z_k)$, the Radon--Nikodym derivative of $\Prob$ with respect to $\mathbb{Q}_\delta$ on $\cF_T^{\mathrm{sig}} \vee \cF_k^Z$ (for $k\delta \le T$) is
\begin{equation}
  L_k^\varepsilon
  := \prod_{j=1}^k \frac{\phi_\delta(Z_j - H_j^\varepsilon)}{\phi_\delta(Z_j)}
  = \prod_{j=1}^k G_j^\varepsilon(\mathcal{X}_j^\varepsilon; Z_j),
  \label{eq:discrete-likelihood-ratio}
\end{equation}
where the single-block relative likelihood multiplier is
\begin{equation}
  G_k^\varepsilon(\mathcal{X}_k^\varepsilon; z)
  := \frac{\phi_\delta(z - H_k^\varepsilon)}{\phi_\delta(z)}
  = \exp\!\left\{z^\top R_\delta^{-1} H_k^\varepsilon
    - \frac{1}{2} (H_k^\varepsilon)^\top R_\delta^{-1} H_k^\varepsilon\right\}.
  \label{eq:full-multiplier}
\end{equation}
Conditional on every signal path, each factor $\phi_\delta(z - H_j^\varepsilon)/\phi_\delta(z)$ integrates to one against $\phi_\delta(z)\dd z$, so $L_k^\varepsilon$ defines a valid change of measure without requiring exponential moment assumptions on $H_k^\varepsilon$.

Formally, for bounded measurable test functions $\Phi:\R^{m+n}\to\R$ and $\varphi:\R^m\to\R$, the full filter and its slow marginal are given by
\begin{align}
  \Pi_k^\varepsilon(\Phi)
  &:= \E^{\Prob}[\Phi(X_{t_k}^\varepsilon, Y_{t_k}^\varepsilon) \mid \cF_k^Z],
  \label{eq:full-posterior}\\
  \pi_k^\varepsilon(\varphi)
  &:= (\ProjX \Pi_k^\varepsilon)(\varphi)
  = \Pi_k^\varepsilon(\varphi \otimes 1)
  = \E^{\Prob}[\varphi(X_{t_k}^\varepsilon) \mid \cF_k^Z],
  \label{eq:slow-marginal}
\end{align}
where $\ProjX:\cP(\R^{m+n}) \to \cP(\R^m)$ denotes the slow marginal projection $(\ProjX \Pi)(\dd x) := \Pi(\dd x \times \R^n)$. Our objective is to approximate $\pi_k^\varepsilon$ over $k\delta \le T$ by an autonomous $m$-dimensional filter, while noting that posterior convergence does not require statistical identifiability of all $m$ slow coordinates from the $d$ observations.

Because $H_k^\varepsilon$ depends on the entire path segment $\mathcal{X}_k^\varepsilon$ over $I_k$, Bayesian updating cannot proceed by pointwise multiplication of a terminal-state likelihood. Instead, for $(x,y)\in\R^{m+n}$, $A \in \cB(\R^{m+n})$, and $z\in\R^d$, define the unnormalized path-integral transition kernel
\begin{equation}
  \mathsf{Q}_\delta^\varepsilon((x,y), A; z)
  := \E_{x,y}\!\left[\mathbf{1}_{\{(X_\delta^\varepsilon, Y_\delta^\varepsilon) \in A\}} G_1^\varepsilon(\mathcal{X}_1^\varepsilon; z)\right].
  \label{eq:rescaled-full-kernel}
\end{equation}
For any joint prior $\Pi \in \cP(\R^{m+n})$, the one-step Bayes update map is
\begin{equation}
  \Phi_{\delta,z}^\varepsilon(\Pi)(A)
  := \frac{\int_{\R^{m+n}} \mathsf{Q}_\delta^\varepsilon((x,y), A; z)\,\Pi(\dd x,\dd y)}
          {\int_{\R^{m+n}} \mathsf{Q}_\delta^\varepsilon((x,y), \R^{m+n}; z)\,\Pi(\dd x,\dd y)}.
  \label{eq:full-joint-one-step-map}
\end{equation}
Completing the square in \eqref{eq:full-multiplier} gives $0 < G_1^\varepsilon(\mathcal{X}_1^\varepsilon; z) \le \exp(\frac{1}{2} z^\top R_\delta^{-1} z)$ for every fixed $z \in \R^d$, so the denominator in \eqref{eq:full-joint-one-step-map} is strictly positive and finite. The exact filter recursion is therefore
\begin{equation}
  \Pi_k^\varepsilon = \Phi_{\delta,Z_k}^\varepsilon(\Pi_{k-1}^\varepsilon),
  \qquad
  \pi_k^\varepsilon = \ProjX \Pi_k^\varepsilon,
  \qquad k \ge 1.
  \label{eq:full-filter-recursion}
\end{equation}

\subsection{Generators, Primitive Assumptions, and Frozen Ergodicity}
\label{subsec:generators-and-ergodicity}

The infinitesimal generator of the coupled diffusion \eqref{eq:slow}--\eqref{eq:fast} decomposes across scales as
\begin{equation}
  \cL^\varepsilon \Phi(x,y)
  = \varepsilon^{-1} \cL_x^F \Phi(x,y) + \cL^S \Phi(x,y),
  \qquad \Phi \in C^2(\R^{m+n}),
  \label{eq:full-generator}
\end{equation}
where the slow operator $\cL^S$ and the frozen fast generator $\cL_x^F$ are
\begin{align}
  \cL^S \Phi(x,y) &:= f(x,y)\cdot \nabla_x \Phi(x,y) + \frac{1}{2} a(x,y):D_x^2\Phi(x,y),
  \label{eq:slow-generator}\\
  \cL_x^F \phi(y) &:= b(x,y)\cdot \nabla_y \phi(y) + \frac{1}{2} A_F(x,y):D_y^2\phi(y),
  \qquad \phi \in C^2(\R^n).
  \label{eq:fast-generator}
\end{align}
Independence of $W^X$ and $W^Y$ eliminates mixed slow--fast second-order derivatives; correlated signal noises would introduce a cross-diffusion term and are not treated here.

We fix a weight exponent $p_* := m + 3$, an integer Sobolev regularity index $s > 4 + m/2$, and a baseline localization threshold $R_{\mathrm{loc}} \ge 1$.
We need the following primitive coefficient and moment hypotheses.

\begin{assumption}
\label{ass:dissipative}
The coefficients of \eqref{eq:slow}--\eqref{eq:fast} satisfy the following conditions:
\begin{enumerate}
\item[\textup{(D1)}] \emph{Local smoothness.}
The maps $f:\R^{m+n}\to\R^m$, $\sigma_X:\R^{m+n}\to\R^{m\times m}$, $b:\R^{m+n}\to\R^n$, and $\tau:\R^{m+n}\to\R^{n\times n}$ are locally Lipschitz, and their mixed partial derivatives $D_x^\alpha D_y^\gamma$ exist and are continuous for all multi-indices $|\alpha| \le s+3$ and $|\gamma| \le 3$.

\item[\textup{(D2)}] \emph{Uniform fast dissipativity and ellipticity.}
There exist constants $\kappa, C_0, \lambda, \Lambda > 0$ such that
\begin{equation}
  \langle y, b(x,y) \rangle \le -\kappa|y|^2 + C_0,
  \qquad
  \lambda I_n \preceq A_F(x,y) \preceq \Lambda I_n,
  \qquad (x,y) \in \R^{m+n}.
  \label{eq:fast-dissipativity}
\end{equation}

\item[\textup{(D3)}] \emph{Slow coercivity and polynomial derivative growth.}
There exist $\kappa_s, C_s > 0$ such that
\[
  \langle x, f(x,y) \rangle \le -\kappa_s|x|^2 + C_s(1+|y|^2),
  \qquad
  \|\sigma_X(x,y)\|_{\HS} \le C_s,
  \qquad (x,y) \in \R^{m+n}.
\]
Furthermore, for every $R \ge R_{\mathrm{loc}}$, there exist $C_R, N_R < \infty$ such that
\[
  \sup_{|x| \le R}
  \sum_{g \in \{f,\sigma_X,b,\tau\}}
  \sum_{\substack{|\alpha| \le s+3 \\ |\gamma| \le 3}}
  |D_x^\alpha D_y^\gamma g(x,y)|
  \le C_R (1+|y|)^{N_R},
  \qquad y \in \R^n,
\]
where $|\cdot|$ denotes any fixed Euclidean component norm for vectors and matrices. By the product rule, analogous polynomial bounds hold for $a = \sigma_X\sigma_X^\top$ and $A_F = \tau\tau^\top$.

\item[\textup{(D4)}] \emph{Initial moments.}
Fix $0 < \Theta < \kappa/\Lambda$ and $q \ge \max\{16, 2p_*\}$. The initial laws satisfy
\[
  \sup_{0 < \varepsilon \le 1}
  \E\!\left[|X_0^\varepsilon|^q + e^{\Theta|Y_0^\varepsilon|^2}\right]
  \le M_0 < \infty.
\]
Neither independence of $(X_0^\varepsilon, Y_0^\varepsilon)$ nor stationarity of the initial fast law is required.

\item[\textup{(D5)}] \emph{Local slow ellipticity.}
For every $R \ge R_{\mathrm{loc}}$, there exists $\underline{a}_R > 0$ such that
\[
  a(x,y) \succeq \underline{a}_R I_m,
  \qquad |x| \le R,\quad y \in \R^n.
\]
\end{enumerate}
\end{assumption}

Assumption~\ref{ass:dissipative} imposes conditions solely on the primitive coefficients; smoothness of the averaged coefficients and parameter regularity of the Poisson correctors are derived rather than assumed (Lemma~\ref{lem:cutoff-poisson}).

\begin{lemma}
\label{lem:signal-moments}
Under Assumption~\ref{ass:dissipative}, the coupled system \eqref{eq:slow}--\eqref{eq:fast} admits a unique nonexplosive strong solution for every $0 < \varepsilon \le 1$. Moreover, for every horizon $T < \infty$, there exists $C_{T,M_0} < \infty$ such that
\begin{equation}
  \sup_{0 < \varepsilon \le 1}
  \E\!\left[\sup_{0 \le t \le T} |X_t^\varepsilon|^q\right]
  + \sup_{0 < \varepsilon \le 1} \sup_{0 \le t \le T}
  \E\bigl[V_\Theta(Y_t^\varepsilon)\bigr]
  \le C_{T,M_0}.
  \label{eq:signal-moment-bound}
\end{equation}
\end{lemma}

\begin{proof}
\emph{Step 1.} We first establish nonexplosion and a uniform Gaussian Lyapunov bound for the fast component.
Because the coefficients are locally Lipschitz by \textup{(D1)}, strong existence and pathwise uniqueness hold up to the explosion time. Using \textup{(D2)}--\textup{(D3)}, the quadratic Lyapunov function $U(x,y) := 1 + |x|^2 + |y|^2$ satisfies $\cL^\varepsilon U(x,y) \le C_\varepsilon U(x,y)$, which rules out explosion by Khasminskii's nonexplosion test \cite[Theorem~3.5]{Khasminskii2012}.

To obtain $\varepsilon$-uniform exponential moments for $Y_t^\varepsilon$, fix any $0 < \eta < \kappa/\Lambda$ and compute
\[
  \nabla_y V_\eta(y) = 2\eta y\, V_\eta(y),
  \qquad
  D_y^2 V_\eta(y) = \bigl(2\eta I_n + 4\eta^2 y y^\top\bigr) V_\eta(y).
\]
Substituting these into \eqref{eq:fast-generator} and applying \textup{(D2)} gives, uniformly in $x \in \R^m$,
\begin{align}
  \cL_x^F V_\eta(y)
  &= \bigl[2\eta \langle y, b(x,y)\rangle + \eta \operatorname{tr}(A_F(x,y)) + 2\eta^2 y^\top A_F(x,y) y\bigr] V_\eta(y) \nonumber\\
  &\le \eta\bigl[-2(\kappa - \eta\Lambda)|y|^2 + 2C_0 + n\Lambda\bigr] V_\eta(y)
  \le -c_\eta V_\eta(y) + C_\eta,
  \label{eq:fast-gaussian-lyapunov}
\end{align}
for positive constants $c_\eta, C_\eta$ independent of $x$ (since $\kappa - \eta\Lambda > 0$, the negative quadratic term dominates outside a compact ball). For $\eta = \Theta$, define the exit times $\theta_k := \inf\{t \ge 0 : |X_t^\varepsilon| + |Y_t^\varepsilon| \ge k\}$. Applying It\^o's formula \cite[Section~3.3]{KaratzasShreve1991} to $e^{c_\Theta t/\varepsilon} V_\Theta(Y_t^\varepsilon)$ up to $t \wedge \theta_k$ and taking expectations via the stopped Dynkin formula \cite[Lemma~3.2]{Khasminskii2012} yields
\[
  \E\bigl[e^{c_\Theta (t \wedge \theta_k)/\varepsilon} V_\Theta(Y_{t \wedge \theta_k}^\varepsilon)\bigr]
  \le \E\bigl[V_\Theta(Y_0^\varepsilon)\bigr]
    + \frac{C_\Theta}{\varepsilon}\E\!\left[\int_0^{t\wedge\theta_k} e^{c_\Theta s/\varepsilon}\dd s\right].
\]
Multiplying by $e^{-c_\Theta t/\varepsilon}$ and letting $k \to \infty$ via Fatou's lemma gives
\[
  \E\bigl[V_\Theta(Y_t^\varepsilon)\bigr]
  \le e^{-c_\Theta t/\varepsilon} \E\bigl[V_\Theta(Y_0^\varepsilon)\bigr]
    + \frac{C_\Theta}{c_\Theta}\bigl(1 - e^{-c_\Theta t/\varepsilon}\bigr)
  \le M_0 + \frac{C_\Theta}{c_\Theta},
\]
uniformly in $t \ge 0$ and $0 < \varepsilon \le 1$.

\emph{Step 2.} We next derive the uniform maximal moment bound for the slow component.
Let $U_q(x) := (1+|x|^2)^{q/2}$. Direct differentiation gives $\nabla_x U_q(x) = q(1+|x|^2)^{(q-2)/2} x$ and $|D_x^2 U_q(x)| \le C_q (1+|x|^2)^{(q-2)/2}$. By \textup{(D3)} and Young's inequality,
\[
  \cL^S U_q(x,y)
  \le q(1+|x|^2)^{(q-2)/2}\bigl(-\kappa_s|x|^2 + C_s(1+|y|^2)\bigr) + C_q C_s^2 (1+|x|^2)^{(q-2)/2}
  \le C_q\bigl(U_q(x) + 1 + |y|^q\bigr).
\]
Applying It\^o's formula to $U_q(X_{t\wedge\theta_k}^\varepsilon)$, taking the supremum over $t \in [0,u]$, and bounding the martingale term via the Burkholder--Davis--Gundy (BDG) inequality \cite[Theorem~3.3.28]{KaratzasShreve1991} together with $\|\sigma_X\|_{\HS} \le C_s$ and Young's inequality, we obtain
\begin{align*}
  \E\!\left[\sup_{0\le t\le u\wedge\theta_k} U_q(X_t^\varepsilon)\right]
  &\le \E[U_q(X_0^\varepsilon)]
    + C_q \int_0^u \E\!\left[\sup_{0\le r\le s\wedge\theta_k} U_q(X_r^\varepsilon)\right]\dd s
    + C_q \int_0^u \E\bigl[1 + |Y_s^\varepsilon|^q\bigr]\dd s \\
  &\quad + \frac{1}{2}\E\!\left[\sup_{0\le t\le u\wedge\theta_k} U_q(X_t^\varepsilon)\right]
    + C_q' \int_0^u \E\bigl[U_q(X_{s\wedge\theta_k}^\varepsilon)\bigr]\dd s.
\end{align*}
Absorbing the factor $\frac{1}{2}$ into the left-hand side, applying Gronwall's inequality on $[0,T]$, and letting $k \to \infty$ by monotone convergence yields
\[
  \E\!\left[\sup_{0 \le t \le T} (1 + |X_t^\varepsilon|^2)^{q/2}\right]
  \le C_T\!\left(1 + \E|X_0^\varepsilon|^q + \int_0^T \E\bigl[1 + |Y_t^\varepsilon|^q\bigr]\dd t\right).
\]
Since $|y|^q \le C_{q,\Theta} V_\Theta(y)$, combining this bound with Step~1 and \textup{(D4)} proves \eqref{eq:signal-moment-bound}.
\end{proof}

Let $(P_t^x)_{t\ge 0}$ denote the Markov semigroup of the frozen fast diffusion with generator $\cL_x^F$. Directly from \textup{(D3)} and the inequality $e^{|c_r||y|} \le C_{\eta,c_r} V_\eta(y)$, every uncentered scalar component $g \in \{f_i, a_{ij}, v_r\}$ satisfies, for every $R \ge R_{\mathrm{loc}}$ and every $\eta > 0$,
\begin{equation}
  \sup_{\substack{|x| \le R, y \in \R^n}}
  \frac{\displaystyle \sum_{\substack{|\alpha| \le s+3, |\gamma| \le 3}} |D_x^\alpha D_y^\gamma g(x,y)|}{V_\eta(y)}
  < \infty.
  \label{eq:centered-forcing-bounds}
\end{equation}
The following proposition establishes uniform geometric ergodicity of the frozen fast process and constructs the Poisson correctors $\chi_\Gamma$ that will neutralize the singular $O(\varepsilon^{-1})$ fast generator term in our multiscale expansions \cite{PardouxVeretennikov2001,PavliotisStuart2008}.

\begin{proposition}[Frozen ergodicity and weighted Poisson resolvents]
\label{prop:fast-ergodicity}
Under Assumption~\ref{ass:dissipative}, for every $x \in \R^m$ the frozen fast diffusion admits a unique invariant probability measure $\mu_x$ (with stationary Fokker--Planck density $p_{\mathrm{inv}}(\cdot\mid x)$), and for every $0 < \eta < \kappa/\Lambda$,
\[
  \sup_{x \in \R^m} \mu_x(V_\eta) < \infty.
\]
Moreover, for every $R \ge R_{\mathrm{loc}}$ and $0 < \eta < \kappa/\Lambda$, there exist constants $C_{R,\eta}, \lambda_{R,\eta} > 0$ such that
\begin{equation}
  \bigl|P_t^x \phi(y) - \mu_x(\phi)\bigr|
  \le C_{R,\eta} e^{-\lambda_{R,\eta} t} V_\eta(y)\|\phi\|_{V_\eta},
  \qquad |x| \le R,\quad y \in \R^n,\quad t \ge 0,
  \label{eq:uniform-mixing}
\end{equation}
where $\|\phi\|_{V_\eta} := \sup_{y\in\R^n} |\phi(y)|/V_\eta(y)$. Writing $\bar g(x) := \int_{\R^n} g(x,y)\,\mu_x(\dd y)$ componentwise and defining the centered forcings
\[
  \Gamma_f := f - \bar f,
  \qquad
  \Gamma_a := a - \bar a,
  \qquad
  \Gamma_v := v - \bar v,
\]
each scalar component $\Gamma$ of $\Gamma_f, \Gamma_a, \Gamma_v$ admits the absolutely convergent Poisson resolvent
\begin{equation}
  \chi_\Gamma(x,y) := \int_0^\infty P_t^x \Gamma(x,\cdot)(y)\dd t.
  \label{eq:poisson-resolvent}
\end{equation}
For each fixed $x \in \R^m$, $\chi_\Gamma(x,\cdot)$ is a $\mu_x$-centered classical solution of $-\cL_x^F \chi_\Gamma = \Gamma$, unique among centered classical solutions bounded by $C_\eta V_\eta(y)$ for every sufficiently small $\eta > 0$.
\end{proposition}

\begin{proof}
\emph{Step 1.} We first verify a Foster--Lyapunov drift condition and establish uniform geometric ergodicity.
Fix $0 < \eta < \kappa/\Lambda$. Because $-2(\kappa - \eta\Lambda)|y|^2 + 2C_0 + n\Lambda \to -\infty$ as $|y|\to\infty$, there exist a radius $K > 0$ and constants $c_\eta', B_\eta > 0$ independent of $x$ such that \eqref{eq:fast-gaussian-lyapunov} implies the Foster--Lyapunov drift condition
\begin{equation}
  \cL_x^F V_\eta(y) \le -c_\eta' V_\eta(y) + B_\eta \mathbf{1}_{B_K}(y),
  \qquad (x,y) \in \R^{m+n}.
  \label{eq:foster-lyapunov}
\end{equation}
Integrating \eqref{eq:fast-gaussian-lyapunov} along the frozen diffusion (as in Lemma~\ref{lem:signal-moments}) yields the discrete-time contraction $P_1^x V_\eta(y) \le e^{-c_\eta} V_\eta(y) + C_\eta/c_\eta$, and likewise $P_t^x V_\eta(y) \le C_\eta' V_\eta(y)$ for all $t \in [0,1]$.

For $|x| \le R$, uniform ellipticity \textup{(D2)} and the local coefficient bounds \textup{(D3)} ensure that the transition density of the frozen diffusion killed upon exiting a sufficiently large ball containing $B_K$ is bounded below by a positive constant independent of $x \in \overline{B}_R$. Thus, every $V_\eta$-sublevel set satisfies a uniform Doeblin minorization condition for $P_1^x$, and Harris's ergodic theorem \cite[Theorems~1.2--1.3]{HairerMattingly2011} yields a unique invariant probability measure $\mu_x$ with geometric convergence at integer times uniformly on $|x| \le R$. Combining this with the uniform bound on $P_t^x V_\eta$ for $t \in [0,1]$ extends the contraction to all $t \ge 0$, proving \eqref{eq:uniform-mixing}. Testing \eqref{eq:fast-gaussian-lyapunov} against $\mu_x$ via bounded truncations of $V_\eta$ and monotone convergence then gives $\sup_{x\in\R^m} \mu_x(V_\eta) \le C_\eta/c_\eta < \infty$.

\emph{Step 2.} We now prove existence, classical regularity, and uniqueness of the Poisson corrector.
For any component $\Gamma$ of $\Gamma_f, \Gamma_a, \Gamma_v$, \eqref{eq:centered-forcing-bounds} and the uniform bound on $\mu_x(V_\eta)$ imply $|\Gamma(x,y)| \le C_{R,\eta} V_\eta(y)$ for every $0 < \eta < \kappa/\Lambda$. By \eqref{eq:uniform-mixing}, the integral \eqref{eq:poisson-resolvent} converges absolutely with $|\chi_\Gamma(x,y)| \le C_{R,\eta}' V_\eta(y)$, and Fubini's theorem gives $\mu_x(\chi_\Gamma(x,\cdot)) = 0$. Weighted Poisson solvability follows from \cite[Theorem~3.2]{GlynnMeyn1996} using Lyapunov function $V_\eta/c_\eta'$ in \eqref{eq:foster-lyapunov}. To verify classical $C^2$ regularity in $y$, consider the discounted resolvent $u_\lambda(y) := \int_0^\infty e^{-\lambda t} P_t^x \Gamma(x,\cdot)(y)\dd t$ ($\lambda > 0$), which solves $(\lambda - \cL_x^F) u_\lambda = \Gamma(x,\cdot)$ and is bounded uniformly in $\lambda \in (0,1]$ in $\|\cdot\|_{V_\eta}$. Applying interior Schauder estimates \cite[Chapter~6]{GilbargTrudinger2001} on compact subsets of $\R^n$ and passing to $\lambda \downarrow 0$ identifies $\chi_\Gamma(x,\cdot)$ as a classical solution of $-\cL_x^F \chi_\Gamma = \Gamma$.

For uniqueness, suppose $h$ is the difference of two $\mu_x$-centered classical solutions bounded by $C_{\eta_0} V_{\eta_0}$ for some $0 < \eta_0 < \kappa/\Lambda$, so that $\cL_x^F h = 0$. Choosing $\eta_1 \in (\eta_0, \kappa/\Lambda)$ and exit times $\theta_k := \inf\{t \ge 0 : |Y_t^x| \ge k\}$, stopped It\^o's formula gives $h(y) = \E_y[h(Y_{t\wedge\theta_k}^x)]$. Since $\sup_{s\ge 0} \E_y[V_{\eta_1}(Y_s^x)] < \infty$ provides uniform integrability, letting $k \to \infty$ yields $P_t^x h = h$ for all $t \ge 0$, and sending $t \to \infty$ in \eqref{eq:uniform-mixing} forces $h = \mu_x(h) = 0$.
\end{proof}

\subsection{Homogenized Signal, Effective Observation Map, and Reduced Filter}
\label{subsec:averaged-signal-and-filter}

Averaging the slow drift and diffusion covariance against the frozen invariant measure $\mu_x$ yields the effective coefficients
\begin{equation}
  \bar f(x) := \int_{\R^n} f(x,y)\,\mu_x(\dd y),
  \qquad
  \bar a(x) := \int_{\R^n} a(x,y)\,\mu_x(\dd y).
  \label{eq:effective-coefficients}
\end{equation}
Note that averaging acts on the full diffusion covariance $a = \sigma_X\sigma_X^\top$ rather than on $\sigma_X$ itself. Letting $\bar\sigma(x) := \bar a(x)^{1/2}$ denote the symmetric positive definite square root, the autonomous $m$-dimensional homogenized signal solves
\begin{equation}
  \dd\bar X_t = \bar f(\bar X_t)\dd t + \bar\sigma(\bar X_t)\dd\bar W_t,
  \qquad \bar X_t \in \R^m,
  \label{eq:averaged-signal}
\end{equation}
where $\bar W$ is an $m$-dimensional standard Brownian motion. The parameter regularity proved in Lemma~\ref{lem:cutoff-poisson}, the local ellipticity \textup{(D5)}, and the coercivity estimate in the proof of Theorem~\ref{thm:filter-rate} guarantee that \eqref{eq:averaged-signal} defines a unique conservative diffusion semigroup $(\bar P_t)_{t\ge 0}$ with backward generator and formal forward adjoint
\[
  \bar\cL \varphi := \bar f \cdot \nabla_x \varphi + \frac{1}{2}\bar a : D_x^2\varphi,
  \qquad
  \bar\cL^* \rho := -\sum_{i=1}^m \partial_{x_i}(\bar f_i \rho)
                    + \frac{1}{2}\sum_{i,j=1}^m \partial_{x_i x_j}(\bar a_{ij}\rho).
\]

Recall the effective observation map $\bar H(x) = (\bar H_1(x),\dots,\bar H_d(x))^\top$ with components $\bar H_r(x) = \log \bar v_r(x)$ defined in \eqref{eq:homogenized-realized-measure-map}. Applying Jensen's inequality to the strictly convex exponential function gives
\[
  \bar H_r(x)
  = \ell_r^\top x + \log \int_{\R^n} e^{c_r^\top y}\,\mu_x(\dd y)
  \ge \ell_r^\top x + c_r^\top \int_{\R^n} y\,\mu_x(\dd y),
\]
with strict inequality whenever $c_r^\top Y$ is nonconstant under $\mu_x$. When interpreting $Y_t^\varepsilon$ as zero-mean fast fluctuations, one may optionally impose the centering convention $\int_{\R^n} y\,\mu_x(\dd y) = 0$ for all $x \in \R^m$, under which $\bar H_r(x) > \ell_r^\top x$ strictly. None of our convergence results require this centering condition; it serves only to interpret the convexity shift relative to the baseline level $\ell_r^\top x$.

The reduced endpoint-observation model is
\begin{equation}
  \bar Z_k = \bar H(\bar X_{t_k}) + R_\delta^{1/2}\xi_k,
  \qquad \xi_k \overset{\mathrm{i.i.d.}}{\sim} \mathcal{N}(0,I_d),
  \label{eq:homogenized-obs-explicit}
\end{equation}
with reduced terminal-state likelihood $\phi_\delta(z - \bar H(x))$. For any prior $\nu \in \cP(\R^m)$ (not necessarily possessing a density), the one-step reduced Bayes update map is
\begin{equation}
  \bar\Phi_{\delta,z}(\nu)(\dd x)
  := \frac{\phi_\delta(z - \bar H(x))\, (\bar P_\delta^* \nu)(\dd x)}
          {(\bar P_\delta^* \nu)\bigl(\phi_\delta(z - \bar H(\cdot))\bigr)},
  \qquad z \in \R^d.
  \label{eq:reduced-one-step-map}
\end{equation}
Initialized from the slow marginal $\bar\pi_0 := \ProjX \Pi_0^\varepsilon$ and driven by the actual observations $(Z_k)_{k\ge 1}$ from \eqref{eq:realized-measure-observation}, the reduced filter evolves via $\bar\pi_k := \bar\Phi_{\delta,Z_k}(\bar\pi_{k-1})$. When densities exist, \eqref{eq:reduced-one-step-map} reduces to the standard predict--update recursion:
\begin{align}
  \bar p_{k\mid k-1} &= \bar P_\delta^* \bar p_{k-1},
  \label{eq:discrete-prediction}\\
  \bar p_k(x) &= \frac{\phi_\delta(Z_k - \bar H(x))\,\bar p_{k\mid k-1}(x)}
                      {\int_{\R^m} \phi_\delta(Z_k - \bar H(u))\,\bar p_{k\mid k-1}(u)\dd u}.
  \label{eq:discrete-realized-update}
\end{align}

Comparing \eqref{eq:full-filter-recursion} and \eqref{eq:reduced-one-step-map} isolates the two approximation steps analyzed in Section~\ref{sec:convergence}: replacing the path-integral likelihood $\phi_\delta(z - H_k^\varepsilon)$ with the convexity-adjusted endpoint likelihood $\phi_\delta(z - \bar H(x))$, and replacing the coupled $(m+n)$-dimensional dynamics \eqref{eq:slow}--\eqref{eq:fast} with the autonomous averaged diffusion \eqref{eq:averaged-signal}.

\section{Convergence Theory for the Multiscale Filter}
\label{sec:convergence}

\subsection{Main Convergence Theorem}
\label{subsec:main-theorem-and-scale-balancing}

\begin{theorem}[Multiscale block-filter convergence and error bound]
\label{thm:filter-rate}
Suppose Assumption~\ref{ass:dissipative} and the Gaussian block observation model \eqref{eq:realized-measure-observation} with covariance scaling \eqref{eq:observation-noise-scaling} hold. Initialize the reduced filter with $\bar\pi_0 = \ProjX \Pi_0^\varepsilon$. Fix $m,n,d$, an integer $s > 4 + m/2$, and a finite horizon $T < \infty$. Then for every localization radius $R \ge R_{\mathrm{loc}}$ and all $0 < \varepsilon \le \delta \le 1$,
\begin{equation}
  \begin{aligned}
    \max_{k\delta \le T} \E^{\Prob}\, d_{\BL}(\pi_k^\varepsilon, \bar\pi_k)
    &\le C_{T,R,M_0}\!\left[
      \sqrt{\frac{\varepsilon}{\delta}} + \sqrt{\delta}
      + \left(\frac{\varepsilon}{\delta}\right)^{1/s}\right]
      + 4\vartheta_T(R) \\
    &\le C_{T,R,M_0}\!\left[\sqrt{\frac{\varepsilon}{\delta}} + \sqrt{\delta}\right]^{1/s}
      + 4\vartheta_T(R),
  \end{aligned}
  \label{eq:upgraded-main-theorem}
\end{equation}
where $\vartheta_T(R) := C_{T,M_0} R^{-q}$. In particular, for every fixed $R \ge R_{\mathrm{loc}}$,
\begin{equation}
  \limsup_{\substack{\varepsilon/\delta \to 0 \\ \delta \to 0}}
  \max_{k\delta \le T} \E^{\Prob}\, d_{\BL}(\pi_k^\varepsilon, \bar\pi_k)
  \le 4\vartheta_T(R).
  \label{eq:iterated-limit-fixed-R}
\end{equation}
Letting $R \to \infty$ establishes finite-horizon convergence of the $m$-dimensional reduced posterior $\bar\pi_k$ to the slow marginal $\pi_k^\varepsilon$ of the $(m+n)$-dimensional full filter in expected bounded-Lipschitz distance.
\end{theorem}

\begin{remark}[Scale balancing and competing block errors]
\label{rem:dual-error-mechanisms}
The block-functional error $r_{\varepsilon,\delta} = \sqrt{\varepsilon/\delta} + \sqrt{\delta}$ reflects the tradeoff between incomplete fast averaging over $I_k$ and slow-trajectory displacement over $I_k$. Balancing these two physical contributions yields the geometric-mean window scaling $\delta \asymp \sqrt{\varepsilon}$, or in physical time,
\begin{equation}
  \Delta^* \asymp \sqrt{\tau_f \tau_s}.
  \label{eq:optimal-block-scale}
\end{equation}
This scaling also minimizes the simplified upper bound $r_{\varepsilon,\delta}^{1/s}$ in \eqref{eq:upgraded-main-theorem}. For the sharper first bound in \eqref{eq:upgraded-main-theorem}, balancing the transition-telescoping term $(\varepsilon/\delta)^{1/s}$ against $\sqrt{\delta}$ gives $\delta \asymp \varepsilon^{2/(s+2)}$ with localized upper-bound order $\varepsilon^{1/(s+2)}$. Both scalings balance analytical upper bounds rather than asserting an exact statistical optimum; the U-shaped window tradeoff associated with \eqref{eq:optimal-block-scale} is examined numerically in Section~\ref{sec:unbalanced}.
\end{remark}

\begin{remark}[Initial nonstationarity, localization, and dimension dependence]
\label{rem:theorem-scope}
Under \textup{(D4)}, the initial fast law need not be stationary; its transient contribution enters the boundary terms of \eqref{eq:backward-corrector-transition} and accumulates to $O(\varepsilon/\delta)$ over $[0,T]$ without assuming that the conditional fast distribution equilibrates under the path-integral likelihood. Because the growth of $C_{T,R,M_0}$ in $R$ is not quantified explicitly, Theorem~\ref{thm:filter-rate} provides a localized rate together with a polynomial exit remainder $4\vartheta_T(R)$. Finally, the Sobolev index $s > 4 + m/2$ depends only on the slow dimension $m$; the fast and observation dimensions $n,d$ affect only the constants and the coefficient hypotheses.
\end{remark}

Direct comparison of the posterior probability measures $\pi_k^\varepsilon$ and $\bar\pi_k$ is obstructed by the random normalizing constants in the denominators of \eqref{eq:full-joint-one-step-map} and \eqref{eq:reduced-one-step-map}. Following the reference-measure approach \cite{BainCrisan2009}, we work with unnormalized filtering measures under the product Gaussian reference law $Q_k := \mathbb{Q}_\delta|_{\cF_k^Z}$, where $(Z_1,\dots,Z_k)$ are i.i.d.\ $\mathcal{N}(0,R_\delta)$ vectors independent of the signal.

For any realization $z_{1:k} \in (\R^d)^k$, we introduce five unnormalized finite measures on $\R^m$:
\begin{enumerate}
  \item[\textup{(1)}] $\rho_k^\varepsilon$: the unnormalized slow marginal of the original full multiscale filter;
  \item[\textup{(2)}] $\rho_k^{\varepsilon,R}$: the unnormalized slow marginal of the localized full system \eqref{eq:cutoff-sde-slow}--\eqref{eq:cutoff-sde-fast};
  \item[\textup{(3)}] $\widehat\rho_k^{\varepsilon,R}$: the unnormalized slow marginal under the intermediate endpoint-observation model $\widehat Z_j = \bar H^R(X_{t_j}^{\varepsilon,R}) + R_\delta^{1/2}\xi_j$;
  \item[\textup{(4)}] $\bar\rho_k^R$: the unnormalized measure of the localized reduced filter;
  \item[\textup{(5)}] $\bar\rho_k$: the unnormalized measure of the original (unlocalized) reduced filter.
\end{enumerate}
By the triangle inequality and $d_{\BL}(\mu,\nu) \le \|\mu-\nu\|_{\TV}$, the expected unnormalized discrepancy decomposes into four distinct errors:
\begin{equation}
  \begin{aligned}
    \E^{Q_k}\, d_{\BL}(\rho_k^\varepsilon, \bar\rho_k)
    &\le \underbrace{\E^{Q_k} \|\rho_k^\varepsilon - \rho_k^{\varepsilon,R}\|_{\TV}}_{\text{Stage 1a: Full localization tail}}
     + \underbrace{\E^{Q_k} \|\rho_k^{\varepsilon,R} - \widehat\rho_k^{\varepsilon,R}\|_{\TV}}_{\text{Stage 2: Observation reduction}} \\
    &\quad + \underbrace{\E^{Q_k}\, d_{\BL}(\widehat\rho_k^{\varepsilon,R}, \bar\rho_k^R)}_{\text{Stage 3: Transition telescoping}}
     + \underbrace{\E^{Q_k} \|\bar\rho_k^R - \bar\rho_k\|_{\TV}}_{\text{Stage 1b: Reduced localization tail}}.
  \end{aligned}
  \label{eq:proof-decomposition-chain}
\end{equation}
The proof follows the decomposition \eqref{eq:proof-decomposition-chain} in three steps: Subsection~\ref{subsec:localization-poisson} constructs the slow-state localization and derives weighted parameter regularity for the Poisson correctors; Subsection~\ref{subsec:observation-reduction} bounds Stage~2 in total variation by $C_{T,R,M_0} r_{\varepsilon,\delta}$; and Subsection~\ref{subsec:transition-telescoping} bounds Stage~3 by $C_{T,R,M_0}(\varepsilon/\delta)^{1/s}$ via weighted Sobolev transition telescoping, controls Stages~1a--1b by $2\vartheta_T(R)$ via pathwise coupling, and synthesizes Theorem~\ref{thm:filter-rate}.

The following elementary Bayes quotient lemma converts unnormalized $Q$-expectations into normalized physical-law posterior bounds without requiring inverse moments of the normalizing constants, and will be used in both Subsection~\ref{subsec:transition-telescoping} and Subsection~\ref{subsec:continuous-time-limit}.

\begin{lemma}
\label{lem:bayes-joint-transfer}
Let $\rho(z)$ and $\rho'(z)$ be positive finite measures on $\R^m$ indexed by $z$, with total masses $\mathcal{Z}(z) := \rho(z)(1) > 0$ and $\mathcal{Z}'(z) := \rho'(z)(1) > 0$. Suppose $\E^Q[\mathcal{Z}] = \E^Q[\mathcal{Z}'] = 1$ under a reference probability measure $Q$, and define the physical observation law $P_Z(\dd z) := \mathcal{Z}(z)\,Q(\dd z)$. Then
\begin{equation}
  \E^{P_Z}\, d_{\BL}\!\left(\frac{\rho}{\mathcal{Z}}, \frac{\rho'}{\mathcal{Z}'}\right)
  \le 2\,\E^Q\, d_{\BL}(\rho, \rho').
  \label{eq:normalization-perturbation}
\end{equation}
The same inequality holds with $d_{\BL}$ replaced by $\|\cdot\|_{\TV}$.
\end{lemma}

\begin{proof}
For any test function $\varphi$ with $\|\varphi\|_{\BL} \le 1$, positivity of $\rho'$ gives $|\rho'(\varphi)/\mathcal{Z}'| \le \|\varphi\|_\infty \le 1$. Writing
\[
  \mathcal{Z}\left|\frac{\rho(\varphi)}{\mathcal{Z}} - \frac{\rho'(\varphi)}{\mathcal{Z}'}\right|
  = \left|(\rho - \rho')(\varphi) + \frac{\rho'(\varphi)}{\mathcal{Z}'}(\mathcal{Z}' - \mathcal{Z})\right|
  \le |(\rho - \rho')(\varphi)| + |(\rho - \rho')(1)|
  \le 2\, d_{\BL}(\rho, \rho'),
\]
taking the supremum over a countable determining class of functions with $\|\varphi\|_{\BL} \le 1$, and integrating against $Q$ proves \eqref{eq:normalization-perturbation}. Replacing $\|\varphi\|_{\BL}\le 1$ by $\|\varphi\|_\infty\le 1$ gives the total-variation analogue.
\end{proof}

\subsection{Slow-State Localization and Regularity of Poisson Correctors}
\label{subsec:localization-poisson}

Fix $R \ge R_{\mathrm{loc}}$ and select a smooth radial truncation map $\psi_R:\R^m \to \R^m$ satisfying $\psi_R(x) = x$ on $\{|x| \le R\}$, $|\psi_R(x)| \le R+1$ on $\R^m$, and bounded derivatives of all positive orders. Explicitly, choose a smooth radial cutoff $\zeta_R:[0,\infty)\to[0,1]$ with $\zeta_R \equiv 1$ on $[0,R]$ and $\zeta_R \equiv 0$ on $[R+1,\infty)$, set $h_R(r) := \int_0^r \zeta_R(u)\dd u$, and define
\[
  \psi_R(x) := \frac{h_R(|x|)}{|x|} x \quad (x \ne 0),
  \qquad \psi_R(0) := 0.
\]
For every coefficient $g \in \{f, \sigma_X, b, \tau, a, A_F, v\}$, define its localized version by $g^R(x,y) := g(\psi_R(x), y)$. The localized slow--fast signal solves
\begin{align}
  \dd X_t^{\varepsilon,R}
  &= f^R(X_t^{\varepsilon,R}, Y_t^{\varepsilon,R})\dd t
   + \sigma_X^R(X_t^{\varepsilon,R}, Y_t^{\varepsilon,R})\dd W_t^X,
  \label{eq:cutoff-sde-slow}\\
  \dd Y_t^{\varepsilon,R}
  &= \varepsilon^{-1} b^R(X_t^{\varepsilon,R}, Y_t^{\varepsilon,R})\dd t
   + \varepsilon^{-1/2} \tau^R(X_t^{\varepsilon,R}, Y_t^{\varepsilon,R})\dd W_t^Y.
  \label{eq:cutoff-sde-fast}
\end{align}
Note that $X_t^{\varepsilon,R}$ remains an unconstrained $\R^m$-valued process; only the slow argument entering the coefficients is truncated. This construction achieves two goals simultaneously: (a) on $\{|x| \le R\}$, the localized and original systems coincide identically up to the exit time $\tau_R^\varepsilon := \inf\{t \ge 0 : |X_t^\varepsilon| \ge R\}$; and (b) globally on $\R^m$, the localized averaged coefficients are bounded with bounded derivatives and uniformly elliptic diffusion, ensuring a well-behaved parabolic flow on weighted Sobolev spaces.

Set $\mu_x^R := \mu_{\psi_R(x)}$, $\bar g^R := \bar g \circ \psi_R$ for $g \in \{f, a, v\}$, and $\bar H^R := \bar H \circ \psi_R$. Define the localized generators
\[
  \cL_x^{F,R} := \cL_{\psi_R(x)}^F,
  \qquad
  \cL^{S,R} := f^R \cdot \nabla_x + \frac{1}{2} a^R : D_x^2,
  \qquad
  \cL^{\varepsilon,R} := \varepsilon^{-1}\cL_x^{F,R} + \cL^{S,R},
\]
and the localized averaged generator $\bar\cL^R := \bar f^R \cdot \nabla_x + \frac{1}{2}\bar a^R : D_x^2$, with corresponding semigroups $P_t^{\varepsilon,R}$ and $\bar P_t^R$. For $r=1,\dots,d$, let $v_r^R(x,y) := \exp(\ell_r^\top \psi_R(x) + c_r^\top y)$ and define the localized block observables
\begin{equation}
  A_{k,r}^{\varepsilon,R} := \frac{1}{\delta}\int_{I_k} v_r^R(X_t^{\varepsilon,R}, Y_t^{\varepsilon,R})\dd t,
  \qquad
  H_{k,r}^{\varepsilon,R} := \log A_{k,r}^{\varepsilon,R}.
  \label{eq:cutoff-realized-measure}
\end{equation}
The localized centered forcings $\Gamma_f^R := f^R - \bar f^R$, $\Gamma_a^R := a^R - \bar a^R$, and $\Gamma_v^R := v^R - \bar v^R$ have corresponding Poisson correctors solving
\begin{align}
  -\cL_x^{F,R}\chi_v^R &= \Gamma_v^R, \qquad \mu_x^R(\chi_v^R) = 0, \label{eq:observation-poisson}\\
  -\cL_x^{F,R}\chi_f^R &= \Gamma_f^R, \qquad \mu_x^R(\chi_f^R) = 0, \label{eq:signal-poisson}\\
  -\cL_x^{F,R}\chi_a^R &= \Gamma_a^R, \qquad \mu_x^R(\chi_a^R) = 0. \label{eq:diffusion-poisson}
\end{align}

\begin{lemma}
\label{lem:cutoff-poisson}
Under Assumption~\ref{ass:dissipative}, let $g$ be any scalar component of $f$, $a$, or $v$. Then the localized Poisson corrector $\chi_g^R(x,y) = \chi_g(\psi_R(x), y)$ is the unique $\mu_x^R$-centered classical solution of
\begin{equation}
  -\cL_x^{F,R}\chi_g^R = g^R - \bar g^R,
  \qquad \mu_x^R(\chi_g^R) = 0,
  \label{eq:cutoff-poisson-general}
\end{equation}
bounded by $C_\eta V_\eta(y)$ for every sufficiently small $\eta > 0$, and satisfies the resolvent formula
\begin{equation}
  \chi_g^R(x,y) = \int_0^\infty P_t^{x,R}\bigl(g^R(x,\cdot) - \bar g^R(x)\bigr)(y)\dd t.
  \label{eq:cutoff-poisson-resolvent}
\end{equation}
Moreover, for every $\eta > 0$, there exists $C_{R,\eta} < \infty$ such that
\begin{equation}
  \sup_{x \in \R^m}
  \sum_{\substack{|\alpha| \le s+3 \\ |\gamma| \le 2}}
  \bigl|D_x^\alpha D_y^\gamma \chi_g^R(x,y)\bigr|
  \le C_{R,\eta} V_\eta(y),
  \qquad y \in \R^n,
  \label{eq:poisson-corrector-bounds}
\end{equation}
and consequently
\begin{equation}
  \bigl|\cL^{S,R}\chi_g^R(x,y)\bigr| \le C_{R,\eta} V_\eta(y),
  \qquad (x,y) \in \R^{m+n}.
  \label{eq:cutoff-slow-action-bounds}
\end{equation}
Finally, the unlocalized averages $\bar f$, $\bar a$, $\bar v$, and $\bar H$ belong to $C_{\mathrm{loc}}^{s+3}(\R^m)$.
\end{lemma}

\begin{proof}
To circumvent differentiating the invariant density $p_{\mathrm{inv}}(\cdot\mid x)$, we establish regularity in four steps: weighted local Schauder bounds in $y$ (Step~1), parameter differentiation of an anchored resolvent $S_x$ (Step~2), invariant-average response and centering restoration (Step~3), and composition with the truncation $\psi_R$ (Step~4). We first work with the untruncated corrector $\chi_g(x,y)$ uniformly on $|x| \le R+1$ (enlarging the ball slightly when taking parameter differences), and then compose with $\psi_R$. Write $L_x := \cL_x^F$ and $k_* := s+3$. By \eqref{eq:centered-forcing-bounds}, for every $\eta > 0$,
\[
  \sup_{|x| \le R+1}
  \sum_{\substack{|\alpha| \le k_* \\ |\gamma| \le 3}}
  \bigl|D_x^\alpha D_y^\gamma g(x,y)\bigr|
  \le C_{R,\eta} V_\eta(y).
\]

\emph{Step 1.} We first derive weighted interior Schauder estimates in the fast variable.
Fix a H\"older exponent $\alpha_0 \in (0,1)$ and define the weighted local H\"older norms
\[
  \|h\|_{j,\alpha_0;\eta} := \sup_{z \in \R^n} e^{-\eta|z|^2}\|h\|_{C^{j,\alpha_0}(B_1(z))},
  \qquad j \in \{0,2\}.
\]
Because the coefficients of $L_x$ grow polynomially in $y$, global H\"older norms need not be bounded. Fix $y_0 \in \R^n$ and choose a small scaling radius $r(y_0) := c_R \langle y_0 \rangle^{-N}$, where $\langle y \rangle := (1+|y|^2)^{1/2}$ and $N \ge 1$ is chosen larger than the polynomial growth degree $N_{R+1}$ in \textup{(D3)}. Under the rescaling $y = y_0 + r(y_0)\tilde y$ for $\tilde y \in B_1(0)$, the rescaled function $\tilde u(\tilde y) := u(y_0 + r(y_0)\tilde y)$ satisfies
\begin{align*}
  \widetilde{L}_{x,y_0}\tilde u(\tilde y)
  &:= \frac{1}{2} A_F(x, y_0 + r(y_0)\tilde y) : D_{\tilde y}^2 \tilde u(\tilde y)
    + r(y_0) b(x, y_0 + r(y_0)\tilde y) \cdot \nabla_{\tilde y}\tilde u(\tilde y) \\
  &= r(y_0)^2 (L_x u)(y_0 + r(y_0)\tilde y).
\end{align*}
Because $A_F$ is uniformly elliptic and bounded by \textup{(D2)}, and the factor $r(y_0)$ neutralizes the polynomial growth of $b$ and of the $y$-derivatives of $(b, A_F)$, the $C^{\alpha_0}(B_1(0))$-norms of the coefficients of $\widetilde{L}_{x,y_0}$ are bounded uniformly in $y_0 \in \R^n$ and $|x| \le R+1$. Interior Schauder estimates \cite[Theorem~6.2]{GilbargTrudinger2001} therefore apply on $B_{1/2}(0) \subset B_1(0)$ with constants independent of $y_0$. Rescaling back to $B_{r(y_0)/2}(y_0)$ introduces at most polynomial factors in $\langle y_0 \rangle$. Since $r(y_0)|y_0| \le c_R$, the weight $V_\eta(y)$ varies by at most a uniform multiplicative constant on $B_{r(y_0)}(y_0)$. Covering $B_1(z)$ by such balls and absorbing the polynomial factors into any strict increment $\eta_1 - \eta_0 > 0$ yields, for all $0 < \eta_0 < \eta_1 < \kappa/\Lambda$,
\begin{equation}
  \|u\|_{2,\alpha_0;\eta_1}
  \le C_{R,\eta_0,\eta_1}\bigl(\|u\|_{V_{\eta_0}} + \|L_x u\|_{0,\alpha_0;\eta_0}\bigr).
  \label{eq:weighted-fast-schauder}
\end{equation}
Applying \eqref{eq:weighted-fast-schauder} first to the discounted resolvents $u_\lambda(y) = \int_0^\infty e^{-\lambda t}(P_t^x g_x(y) - \mu_x(g_x))\dd t$ (which satisfy $(\lambda - L_x)u_\lambda = g_x - \mu_x(g_x)$) and passing to $\lambda \downarrow 0$ shows that the centered resolvent $\chi_g(x,\cdot)$ is bounded in $\|\cdot\|_{2,\alpha_0;\eta_1}$ uniformly in $|x| \le R+1$.

\emph{Step 2.} We next introduce an anchored resolvent to differentiate with respect to the slow parameter inductively.
To differentiate with respect to $x$ without differentiating the invariant density $p_{\mathrm{inv}}(\cdot\mid x)$, we adapt the transfer-identity method of R\"ockner and Xie \cite{RocknerXie2021} using an anchored resolvent at $y=0$. For any forcing $h$ with finite $\|\cdot\|_{0,\alpha_0;\eta_0}$-norm, define
\begin{equation}
  S_x h(y) := \int_0^\infty \bigl(P_t^x h(y) - P_t^x h(0)\bigr)\dd t,
  \qquad
  u_x := S_x g_x,
  \label{eq:anchored-poisson-resolvent}
\end{equation}
where $g_x(y) := g(x,y)$. Then $u_x(0) = 0$ and $-L_x u_x = g_x - \mu_x(g_x)$. By Step~1, $S_x$ maps the weighted space $\|\cdot\|_{0,\alpha_0;\eta_0}$ continuously into $\|\cdot\|_{2,\alpha_0;\eta_1}$ for any strict increment $0 < \eta_0 < \eta_1 < \kappa/\Lambda$, and multiplication by polynomially growing coefficient derivatives is likewise bounded between weighted spaces with any positive increment in $\eta$.

For two nearby parameters $x, z$ with $|x|, |z| \le R+1$, apply $-L_x$ to $u_z - u_x$:
\[
  -L_x(u_z - u_x)
  = -L_z u_z + L_x u_x + (L_z - L_x)u_z
  = \bigl[g_z - g_x + (L_z - L_x)u_z\bigr] - \bigl(\mu_z(g_z) - \mu_x(g_x)\bigr).
\]
Integrating $-L_x(u_z - u_x)$ against $\mu_x$ (justified by stopped It\^o's formula using a slightly larger Gaussian weight below $\kappa/\Lambda$ for uniform integrability) gives $\mu_x[g_z - g_x + (L_z - L_x)u_z] = \mu_z(g_z) - \mu_x(g_x)$. Since both $u_z - u_x$ and $S_x[g_z - g_x + (L_z - L_x)u_z]$ vanish at $y=0$ and solve the same Poisson equation in the $V_{\eta_1}$-growth class, uniqueness (Proposition~\ref{prop:fast-ergodicity}) yields the exact resolvent difference identity
\begin{equation}
  u_z - u_x = S_x\bigl[g_z - g_x + (L_z - L_x)u_z\bigr].
  \label{eq:anchored-poisson-difference}
\end{equation}

We now choose finitely many strictly increasing Gaussian exponents $0 < \eta_0 < \eta_1 < \dots < \eta_K < \eta$ (for any prescribed target $\eta \in (0,\kappa/\Lambda)$) and strictly decreasing H\"older exponents $\alpha_0 > \alpha_1 > \dots > \alpha_K > 0$. By \eqref{eq:anchored-poisson-difference} and the uniform operator bounds on $S_x$, $x \mapsto u_x$ is continuous in $\|\cdot\|_{2,\alpha_1;\eta_1}$. Dividing \eqref{eq:anchored-poisson-difference} by $z_i - x_i$ (with $z = x + \theta e_i$, $\theta \to 0$) and noting that on bounded sets in $y$ the difference quotients converge locally while outside bounded sets the positive increment in $\eta$ controls the tails uniformly, we obtain the first and second parameter-derivative identities:
\begin{align}
  \partial_{x_i} u_x &= S_x F_i, & F_i &:= \partial_{x_i} g_x + (\partial_{x_i} L_x) u_x,
  \label{eq:poisson-first-parameter-derivative}\\
  \partial_{x_i x_j} u_x &= S_x F_{ij}, & F_{ij} &:= \partial_{x_i x_j} g_x + (\partial_{x_i x_j} L_x) u_x + (\partial_{x_i} L_x)\partial_{x_j} u_x + (\partial_{x_j} L_x)\partial_{x_i} u_x.
  \label{eq:poisson-second-parameter-derivative}
\end{align}
Importantly, because $\partial_{x_i} L_x$ differentiates $u_x$ at most twice in $y$, the right-hand sides remain in weighted $C^{\alpha_j}$ spaces without losing fast-variable regularity at successive parameter-differentiation steps. By induction for every $1 \le |\alpha| \le k_*$,
\begin{equation}
  D_x^\alpha u_x
  = S_x\!\left[D_x^\alpha g_x + \sum_{0 < \nu \le \alpha} \binom{\alpha}{\nu} (D_x^\nu L_x) D_x^{\alpha-\nu} u_x\right],
  \label{eq:poisson-parameter-induction}
\end{equation}
with $\|D_x^\alpha u_x\|_{2,\alpha_{|\alpha|};\eta_{|\alpha|}} \le C_{R,\eta}$ uniformly on $|x| \le R+1$.

\emph{Step 3.} We now establish parameter regularity of the invariant averages and restore centering.
From the uniform Lyapunov bound $\sup_x \mu_x(V_\eta) < \infty$, continuity of coefficients, and uniqueness of $\mu_x$, the map $x \mapsto \mu_x$ is weakly continuous. For any parameter-dependent forcing $h_x$ in the above weighted H\"older class, integrating $-L_x(S_x h_x) = h_x - \mu_x(h_x)$ against $\mu_z$ and using $\mu_z(L_z S_x h_x) = 0$ yields the invariant-measure response identity
\begin{equation}
  \mu_z(h_z) - \mu_x(h_x)
  = \mu_z\bigl[h_z - h_x + (L_z - L_x)S_x h_x\bigr].
  \label{eq:invariant-parameter-response}
\end{equation}
Alternatively, evaluating $-L_x u_x = g_x - \bar g(x)$ at the anchor $y=0$ gives the pointwise identity $\bar g(x) = g(x,0) + (L_x u_x)(0)$, which together with Step~2 immediately implies $\bar g \in C_{\mathrm{loc}}^{k_*}(\R^m)$, with $\partial_{x_i}\bar g = \mu_x(F_i)$ and $\partial_{x_i x_j}\bar g = \mu_x(F_{ij})$.

To restore $\mu_x$-centering, observe that $\chi_g(x,y) = u_x(y) - c(x)$ where $c(x) := \mu_x(u_x)$. Applying \eqref{eq:invariant-parameter-response} and \eqref{eq:poisson-parameter-induction} with the forcing $h_x := u_x$ (letting $w_x := S_x u_x$) yields $c \in C_{\mathrm{loc}}^{k_*}(\R^m)$, with
\[
  \partial_{x_i} w_x = S_x\bigl(\partial_{x_i} u_x + (\partial_{x_i} L_x) w_x\bigr),
  \qquad
  \partial_{x_i} c(x) = \mu_x\bigl(\partial_{x_i} u_x + (\partial_{x_i} L_x) w_x\bigr),
\]
and higher derivatives obtained inductively. Thus $\chi_g = u - c$ satisfies the mixed derivative bounds uniformly on $|x| \le R+1$ for every arbitrarily small $\eta > 0$.

\emph{Step 4.} Finally, we compose with the radial truncation $\psi_R$ and bound the slow-generator action.
Composing $\chi_g$ with $\psi_R$ preserves $\mu_x^R$-centering and yields \eqref{eq:cutoff-poisson-resolvent}. Because $\psi_R$ has bounded derivatives of all positive orders, the Fa\`a di Bruno formula gives \eqref{eq:poisson-corrector-bounds}. For \eqref{eq:cutoff-slow-action-bounds}, apply \eqref{eq:poisson-corrector-bounds} with exponent $\eta/2$ and absorb the polynomial growth of $f^R(x,y)$ into $V_{\eta/2}(y)$ (while $a^R$ is bounded). Finally, since $\bar f, \bar a, \bar v \in C_{\mathrm{loc}}^{k_*}(\R^m)$ and each $\bar v_r > 0$, we conclude $\bar H_r = \log \bar v_r \in C_{\mathrm{loc}}^{k_*}(\R^m)$.
\end{proof}

We now introduce the localized observation model
\begin{equation}
  Z_k = H_k^{\varepsilon,R} + R_\delta^{1/2}\xi_k,
  \qquad \xi_k \overset{\mathrm{i.i.d.}}{\sim} \mathcal{N}(0,I_d),
  \label{eq:cutoff-block-observation}
\end{equation}
with joint probability law $\Prob^R$ and single-block relative likelihood
\begin{equation}
  G_k^{\varepsilon,R}(\mathcal{X}_k^{\varepsilon,R}; z)
  := \exp\!\left\{z^\top R_\delta^{-1} H_k^{\varepsilon,R}
    - \frac{1}{2}(H_k^{\varepsilon,R})^\top R_\delta^{-1} H_k^{\varepsilon,R}\right\}.
  \label{eq:cutoff-multiplier}
\end{equation}
For the localized reduced model, define the endpoint relative likelihood
\begin{equation}
  G_k^R(x; z)
  := \frac{\phi_\delta(z - \bar H^R(x))}{\phi_\delta(z)}
  = \exp\!\left\{z^\top R_\delta^{-1}\bar H^R(x)
    - \frac{1}{2}\bar H^R(x)^\top R_\delta^{-1}\bar H^R(x)\right\},
  \label{eq:reduced-cutoff-likelihood}
\end{equation}
together with the unnormalized backward and forward one-step propagators $\bar Q_{k,z}^R u := \bar P_\delta^R(u\, G_k^R(\cdot; z))$ and $\bar Q_k^{R,*}\rho := G_k^R(\cdot; Z_k)\,\bar P_\delta^{R,*}\rho$. The corresponding normalized full and reduced posteriors are denoted by $\Pi_k^{\varepsilon,R}$, $\pi_k^{\varepsilon,R} := \ProjX \Pi_k^{\varepsilon,R}$, and $\bar\pi_k^R$.

\begin{lemma}
\label{lem:cutoff-experiment}
Under Assumption~\ref{ass:dissipative}, the localized signal \eqref{eq:cutoff-sde-slow}--\eqref{eq:cutoff-sde-fast} admits a unique global strong solution with frozen invariant measure $\mu_x^R$. Moreover, $\bar f^R, \bar a^R, \bar H^R \in C_b^{s+3}(\R^m)$, and there exist constants $\underline{a}_R, \overline{a}_R, M_R > 0$ such that
\[
  \underline{a}_R I_m \preceq \bar a^R(x) \preceq \overline{a}_R I_m,
  \qquad |\bar H^R(x)| \le M_R,
  \qquad x \in \R^m.
\]
Consequently, $\bar\cL^R$ generates a conservative diffusion semigroup $\bar P_t^R$, and for every $z \in \R^d$ both likelihood multipliers \eqref{eq:cutoff-multiplier}--\eqref{eq:reduced-cutoff-likelihood} are strictly positive and bounded above by $\exp(\frac{1}{2} z^\top R_\delta^{-1} z)$, so the localized Bayes update maps are well defined.
\end{lemma}

\begin{proof}
Because $\psi_R$ is smooth and bounded, the localized coefficients remain locally Lipschitz, and $f^R(x,y)$ grows at most polynomially in $y$ uniformly in $x \in \R^m$. Choosing $\eta \in (0,\kappa/\Lambda)$ sufficiently small, the Lyapunov function $1 + |x|^2 + V_\eta(y)$ has $\cL^{\varepsilon,R}$-generator bounded above by $C_{R,\varepsilon}(1 + |x|^2 + V_\eta(y))$, ruling out explosion.

By Lemma~\ref{lem:cutoff-poisson} and the boundedness of $\psi_R$ and its derivatives, $\bar f^R, \bar a^R, \bar H^R \in C_b^{s+3}(\R^m)$. Averaging \textup{(D5)} and $\|\sigma_X\|_{\HS} \le C_s$ against $\mu_x^R$ gives $\underline{a}_R I_m \preceq \bar a^R(x) \preceq \overline{a}_R I_m$. Since the symmetric square-root map is smooth on uniformly positive definite matrices, $(\bar a^R)^{1/2}$ is bounded and globally Lipschitz, so $\bar\cL^R$ generates a conservative semigroup. Uniform invariant exponential moments bound $\bar v_r^R(x)$ from above, while Jensen's inequality gives
\[
  \bar v_r^R(x) \ge \exp\!\bigl(\ell_r^\top \psi_R(x) + c_r^\top \mu_x^R(Y)\bigr) \ge c_R > 0,
\]
yielding $|\bar H^R(x)| \le M_R$. Finally, completing the square in the Gaussian exponents gives the upper bound $\exp(\frac{1}{2} z^\top R_\delta^{-1} z)$, and path continuity guarantees $0 < A_{k,r}^{\varepsilon,R} < \infty$ almost surely.
\end{proof}

\begin{lemma}
\label{lem:localized-block-moments}
For every fixed $R \ge R_{\mathrm{loc}}$ and $T < \infty$,
\begin{equation}
  \sup_{0 < \varepsilon \le 1} \sup_{0 \le t \le T} \E\bigl[V_\Theta(Y_t^{\varepsilon,R})\bigr]
  + \sup_{0 < \varepsilon \le 1} \E\!\left[\sup_{0 \le t \le T} |X_t^{\varepsilon,R}|^q\right]
  \le C_{R,T,M_0}.
  \label{eq:cutoff-signal-moments}
\end{equation}
Moreover, for every $p \ge 1$ and $r \in \{1,\dots,d\}$,
\begin{equation}
  \sup_{\substack{0 < \varepsilon \le \delta \le 1 \\ k\delta \le T}}
  \E\!\left[(A_{k,r}^{\varepsilon,R})^p + (A_{k,r}^{\varepsilon,R})^{-p}
    + \bigl(\bar v_r^R(X_{t_k}^{\varepsilon,R})\bigr)^p
    + \bigl(\bar v_r^R(X_{t_k}^{\varepsilon,R})\bigr)^{-p}\right]
  \le C_{p,R,T,M_0}.
  \label{eq:localized-block-moments}
\end{equation}
\end{lemma}

\begin{proof}
The fast Lyapunov inequality \eqref{eq:fast-gaussian-lyapunov} is unaffected by slow-argument truncation, so the first bound in \eqref{eq:cutoff-signal-moments} follows exactly as in Lemma~\ref{lem:signal-moments}. For $X_t^{\varepsilon,R}$, using $|f^R(x,y)| \le C_R(1+|y|^{N_R})$, boundedness of $\sigma_X^R$, Jensen's inequality, and the BDG inequality gives
\[
  \E\!\left[\sup_{0 \le t \le T} |X_t^{\varepsilon,R}|^q\right]
  \le C_{R,T,q}\!\left(1 + \E|X_0^\varepsilon|^q + \int_0^T \E\bigl[1 + |Y_t^{\varepsilon,R}|^{q N_R}\bigr]\dd t\right),
\]
which is finite by the uniform $V_\Theta$-moment bound and \textup{(D4)}.

For \eqref{eq:localized-block-moments}, note that both $u \mapsto u^p$ and $u \mapsto u^{-p}$ are convex on $(0,\infty)$ for $p \ge 1$. Applying Jensen's inequality to the time average $A_{k,r}^{\varepsilon,R} = \delta^{-1}\int_{I_k} \exp(\ell_r^\top \psi_R(X_t^{\varepsilon,R}) + c_r^\top Y_t^{\varepsilon,R})\dd t$ and using $|\psi_R(x)| \le R+1$ yields
\[
  (A_{k,r}^{\varepsilon,R})^{\pm p}
  \le \frac{1}{\delta}\int_{I_k} \exp\!\bigl(\pm p \ell_r^\top \psi_R(X_t^{\varepsilon,R}) \pm p c_r^\top Y_t^{\varepsilon,R}\bigr)\dd t
  \le \frac{C_{p,R}}{\delta}\int_{I_k} e^{p|c_r||Y_t^{\varepsilon,R}|}\dd t.
\]
Completing the square via $p|c_r||y| \le \Theta|y|^2 + p^2|c_r|^2/(4\Theta)$ bounds both expectations by $C_{p,R,\Theta}\sup_{t\le T}\E[V_\Theta(Y_t^{\varepsilon,R})] \le C_{p,R,T,M_0}$. The remaining two terms in \eqref{eq:localized-block-moments} are deterministically bounded by Lemma~\ref{lem:cutoff-experiment}.
\end{proof}

To accommodate all eighth moments needed in the block-observation comparison within our $V_\Theta$ moment budget, we fix the common Gaussian weight envelope
\begin{equation}
  W(y) := V_{\Theta/16}(y),
  \qquad\text{so that}\qquad
  W(y)^4 = V_{\Theta/4}(y),
  \quad
  W(y)^8 = V_{\Theta/2}(y) \le V_\Theta(y).
  \label{eq:common-corrector-envelope}
\end{equation}
Applying Lemma~\ref{lem:cutoff-poisson} with $\eta = \Theta/32$ and absorbing polynomially growing coefficients into another factor $V_{\Theta/32}$ yields
\begin{equation}
  |\chi_v^R| + \|D_x\chi_v^R \sigma_X^R\|_{\HS} + \|D_y\chi_v^R \tau^R\|_{\HS} + |\cL^{S,R}\chi_v^R|
  \le C_R W.
  \label{eq:observation-corrector-envelope}
\end{equation}

\subsection{Block-to-Endpoint Observation Reduction}
\label{subsec:observation-reduction}

We now quantify the error incurred by replacing the path-dependent block observable $H_k^{\varepsilon,R}$ with the endpoint observable $\bar H^R(X_{t_k}^{\varepsilon,R})$.

\begin{lemma}
\label{lem:fast-block}
For every $p \in [2,8]$, $0 < \varepsilon \le \delta \le 1$, and $k\delta \le T$,
\begin{equation}
  \left\|\frac{1}{\delta}\int_{I_k} \Gamma_v^R(X_t^{\varepsilon,R}, Y_t^{\varepsilon,R})\dd t\right\|_{L^p}
  \le C_{p,R,T,M_0}\!\left(\sqrt{\frac{\varepsilon}{\delta}} + \frac{\varepsilon}{\delta} + \varepsilon\right)
  \le C_{p,R,T,M_0}'\sqrt{\frac{\varepsilon}{\delta}}.
  \label{eq:block-mixing}
\end{equation}
\end{lemma}

\begin{proof}
Applying It\^o's formula componentwise to $\chi_v^R(X_t^{\varepsilon,R}, Y_t^{\varepsilon,R})$ on $I_k = (t_{k-1}, t_k]$ and using the Poisson equation $-\cL_x^{F,R}\chi_v^R = \Gamma_v^R$ together with the independence of $W^X$ and $W^Y$ (which eliminates any mixed second-order term), we obtain
\begin{align}
  \int_{I_k} \Gamma_v^R(X_t^{\varepsilon,R}, Y_t^{\varepsilon,R})\dd t
  &= -\varepsilon\bigl[\chi_v^R(X_t^{\varepsilon,R}, Y_t^{\varepsilon,R})\bigr]_{t_{k-1}}^{t_k}
   + \varepsilon \int_{I_k} \cL^{S,R}\chi_v^R(X_t^{\varepsilon,R}, Y_t^{\varepsilon,R})\dd t \nonumber\\
  &\quad + \varepsilon \int_{I_k} D_x\chi_v^R\,\sigma_X^R \dd W_t^X
   + \sqrt{\varepsilon} \int_{I_k} D_y\chi_v^R\,\tau^R \dd W_t^Y.
  \label{eq:corrector-identity}
\end{align}
(Strictly, \eqref{eq:corrector-identity} is derived up to stopping times and passed to the limit via \eqref{eq:observation-corrector-envelope} and Lemma~\ref{lem:localized-block-moments}.) For any $p \in [2,8]$, we have $W(y)^p \le V_\Theta(y)$, so $\sup_{t\le T}\|W(Y_t^{\varepsilon,R})\|_{L^p} \le C_{p,R,T,M_0}$. Bounding the four terms on the right-hand side of \eqref{eq:corrector-identity} in $L^p$ via \eqref{eq:observation-corrector-envelope}, Minkowski's integral inequality, and the BDG inequality gives:
\begin{align*}
  \varepsilon\bigl\|\bigl[\chi_v^R\bigr]_{t_{k-1}}^{t_k}\bigr\|_{L^p} &\le C\varepsilon, &
  \varepsilon\left\|\int_{I_k} \cL^{S,R}\chi_v^R \dd t\right\|_{L^p} &\le C\varepsilon\delta, \\
  \varepsilon\left\|\int_{I_k} D_x\chi_v^R\sigma_X^R \dd W_t^X\right\|_{L^p} &\le C\varepsilon\sqrt{\delta}, &
  \sqrt{\varepsilon}\left\|\int_{I_k} D_y\chi_v^R\tau^R \dd W_t^Y\right\|_{L^p} &\le C\sqrt{\varepsilon\delta},
\end{align*}
where we used $\E[(\int_{I_k} W^2\dd t)^{p/2}] \le \delta^{p/2-1}\int_{I_k}\E[W^p]\dd t \le C\delta^{p/2}$. Dividing by $\delta$ and noting that $\varepsilon \le \delta \le 1$ implies $\varepsilon/\delta \le \sqrt{\varepsilon/\delta}$ and $\varepsilon \le \sqrt{\varepsilon/\delta}$ yields \eqref{eq:block-mixing}.
\end{proof}

\begin{lemma}
\label{lem:slow-block}
For every $p \in [2,8]$, $0 < \varepsilon \le \delta \le 1$, and $k\delta \le T$,
\begin{equation}
  \left\|\frac{1}{\delta}\int_{I_k} \bigl\{\bar v^R(X_t^{\varepsilon,R}) - \bar v^R(X_{t_{k-1}}^{\varepsilon,R})\bigr\}\dd t\right\|_{L^p}
  \le C_{p,R,T,M_0}\sqrt{\delta}.
  \label{eq:slow-block-bound}
\end{equation}
\end{lemma}

\begin{proof}
Using $|f^R(x,y)| \le C_R(1+|y|^{N_R})$, boundedness of $\sigma_X^R$, Lemma~\ref{lem:localized-block-moments}, and the BDG inequality, we have for all $0 \le u \le t \le T$ with $t-u \le 1$,
\[
  \|X_t^{\varepsilon,R} - X_u^{\varepsilon,R}\|_{L^p}
  \le C_{p,R,T,M_0}\bigl((t-u) + \sqrt{t-u}\bigr)
  \le C_{p,R,T,M_0}'\sqrt{t-u}.
\]
Because $\bar v^R$ is globally Lipschitz by Lemma~\ref{lem:cutoff-poisson}, Minkowski's integral inequality yields \eqref{eq:slow-block-bound}. The identical increment bound holds with $\bar H^R$ in place of $\bar v^R$.
\end{proof}

\begin{lemma}
\label{lem:realized-measure-reduction}
For either block endpoint $i \in \{k-1, k\}$, we have
\begin{equation}
  \bigl\|H_k^{\varepsilon,R} - \bar H^R(X_{t_i}^{\varepsilon,R})\bigr\|_{L^4}
  \le C_{R,T,M_0}\!\left(\sqrt{\frac{\varepsilon}{\delta}} + \sqrt{\delta}\right)
  = C_{R,T,M_0}\, r_{\varepsilon,\delta}.
  \label{eq:realized-measure-homogenization}
\end{equation}
\end{lemma}

\begin{proof}
Decompose the difference of block-averaged and frozen-endpoint variances as
\[
  A_k^{\varepsilon,R} - \bar v^R(X_{t_{k-1}}^{\varepsilon,R})
  = \frac{1}{\delta}\int_{I_k} \Gamma_v^R(X_t^{\varepsilon,R}, Y_t^{\varepsilon,R})\dd t
  + \frac{1}{\delta}\int_{I_k} \bigl\{\bar v^R(X_t^{\varepsilon,R}) - \bar v^R(X_{t_{k-1}}^{\varepsilon,R})\bigr\}\dd t.
\]
By Lemmas~\ref{lem:fast-block} and \ref{lem:slow-block} with $p=8$, the left-hand side has $L^8$ norm at most $C_{R,T,M_0} r_{\varepsilon,\delta}$. For each component $r \in \{1,\dots,d\}$, the elementary inequality $|\log u_1 - \log u_2| \le |u_1 - u_2|(u_1^{-1} + u_2^{-1})$ ($u_1, u_2 > 0$) and the Cauchy--Schwarz inequality give
\begin{align*}
  \bigl\|\log A_{k,r}^{\varepsilon,R} - \log \bar v_r^R(X_{t_{k-1}}^{\varepsilon,R})\bigr\|_{L^4}
  &\le \bigl\|A_{k,r}^{\varepsilon,R} - \bar v_r^R(X_{t_{k-1}}^{\varepsilon,R})\bigr\|_{L^8}
       \Bigl(\bigl\|(A_{k,r}^{\varepsilon,R})^{-1}\bigr\|_{L^8}
       + \bigl\|(\bar v_r^R(X_{t_{k-1}}^{\varepsilon,R}))^{-1}\bigr\|_{L^8}\Bigr).
\end{align*}
The inverse $L^8$ moments are uniformly bounded by Lemma~\ref{lem:localized-block-moments}. Summing over $r=1,\dots,d$ yields \eqref{eq:realized-measure-homogenization} for $i=k-1$, and the $L^4$ increment bound for $\bar H^R(X^{\varepsilon,R})$ from Lemma~\ref{lem:slow-block} gives the right-endpoint case $i=k$.
\end{proof}

On the same localized signal $(X^{\varepsilon,R}, Y^{\varepsilon,R})$, recall the intermediate endpoint-observation model $\widehat Z_j = \bar H^R(X_{t_j}^{\varepsilon,R}) + R_\delta^{1/2}\xi_j$ ($1 \le j \le k$) with unnormalized slow marginal $\widehat\rho_k^{\varepsilon,R}$.

\begin{lemma}
\label{lem:block-endpoint-experiment}
Under Assumption~\ref{ass:dissipative}, for all $0 < \varepsilon \le \delta \le 1$,  we have
\begin{equation}
  \max_{k\delta \le T} \E^{\mathbb{Q}_\delta}
  \bigl\|\rho_k^{\varepsilon,R} - \widehat\rho_k^{\varepsilon,R}\bigr\|_{\TV}
  \le C_{T,R,M_0}\, r_{\varepsilon,\delta}.
  \label{eq:block-endpoint-total-variation}
\end{equation}
\end{lemma}

\begin{proof}
Let $L_k^{\varepsilon,R} := \prod_{j=1}^k G_j^{\varepsilon,R}(\mathcal{X}_j^{\varepsilon,R}; z_j)$ and $\widehat L_k^{\varepsilon,R} := \prod_{j=1}^k G_j^R(X_{t_j}^{\varepsilon,R}; z_j)$ denote the multiblock relative likelihoods under $\mathbb{Q}_\delta$. For any $\|\varphi\|_\infty \le 1$,
\[
  \bigl|(\rho_k^{\varepsilon,R} - \widehat\rho_k^{\varepsilon,R})(\varphi)(z_{1:k})\bigr|
  = \bigl|\E_{\mathrm{sig}}\bigl[\varphi(X_{t_k}^{\varepsilon,R})(L_k^{\varepsilon,R} - \widehat L_k^{\varepsilon,R})\bigr]\bigr|
  \le \E_{\mathrm{sig}}\bigl|L_k^{\varepsilon,R} - \widehat L_k^{\varepsilon,R}\bigr|.
\]
Taking the supremum over $\|\varphi\|_\infty \le 1$ and applying Tonelli's theorem yields
\[
  \E^{\mathbb{Q}_\delta}\bigl\|\rho_k^{\varepsilon,R} - \widehat\rho_k^{\varepsilon,R}\bigr\|_{\TV}
  \le \E_{\mathrm{sig}}\E^{\mathbb{Q}_\delta}\bigl|L_k^{\varepsilon,R} - \widehat L_k^{\varepsilon,R}\bigr|.
\]
Conditional on the signal path, the two observation vectors on $\R^{kd}$ are Gaussian with identical block-diagonal covariance $I_k \otimes R_\delta$ and mean difference $(\Delta_1,\dots,\Delta_k)$, where $\Delta_j := H_j^{\varepsilon,R} - \bar H^R(X_{t_j}^{\varepsilon,R})$. The Hellinger affinity between two Gaussian distributions $\mathcal{N}(m_1,\Sigma)$ and $\mathcal{N}(m_2,\Sigma)$ equals $\exp(-\frac{1}{8}(m_1-m_2)^\top \Sigma^{-1}(m_1-m_2))$. Putting $\mathcal{D}_k := \sum_{j=1}^k \Delta_j^\top R_\delta^{-1} \Delta_j$, the Le Cam inequality relating total variation to Hellinger affinity \cite[Lemma~2.3]{Tsybakov2009} yields
\[
  \E^{\mathbb{Q}_\delta}\bigl|L_k^{\varepsilon,R} - \widehat L_k^{\varepsilon,R}\bigr|
  \le 2\sqrt{1 - e^{-\mathcal{D}_k/4}}
  \le \sqrt{\mathcal{D}_k},
\]
where we used $1 - e^{-u} \le u$ for $u \ge 0$. Since $R_\delta^{-1} \preceq (\delta/r_0)I_d$ by \eqref{eq:observation-noise-scaling}, Cauchy--Schwarz gives
\[
  \E^{\mathbb{Q}_\delta}\bigl\|\rho_k^{\varepsilon,R} - \widehat\rho_k^{\varepsilon,R}\bigr\|_{\TV}
  \le \E_{\mathrm{sig}}\bigl[\sqrt{\mathcal{D}_k}\bigr]
  \le \left(\frac{\delta}{r_0}\sum_{j=1}^k \E|\Delta_j|^2\right)^{1/2}.
\]
Applying Lemma~\ref{lem:realized-measure-reduction} to each block and using $k\delta \le T$ yields \eqref{eq:block-endpoint-total-variation}.
\end{proof}

\begin{remark}[Stability of multiblock observation accumulation]
Because $\|R_\delta^{-1}\|_{\op} \le \delta/r_0$, each block's Fisher information scales as $O(\delta)$, which precisely cancels the $O(1/\delta)$ growth in the number of blocks $k \le T/\delta$. Consequently, $\E[\mathcal{D}_k] \le (T/r_0)\max_{j\le k}\E|\Delta_j|^2 \le C_{T,R,M_0} r_{\varepsilon,\delta}^2$, so the cumulative multiblock observation error scales as $O(r_{\varepsilon,\delta})$ uniformly over $[0,T]$ rather than diverging as $\delta \downarrow 0$.
\end{remark}

\subsection{Proof of Theorem~\ref{thm:filter-rate}}
\label{subsec:transition-telescoping}

To compare the endpoint-observation filter $\widehat\rho_k^{\varepsilon,R}$ with the localized reduced filter $\bar\rho_k^R$, we work in polynomially weighted Sobolev spaces on $\R^m$. Recall $p_* = m + 3$ and the integer $s > 4 + m/2$, and define the integrable polynomial weight
\[
  w(x) := (1 + |x|^2)^{-p_*/2}, \qquad x \in \R^m.
\]
Let $H_w^s(\R^m) := H^s(\R^m, w\dd x)$ be the Hilbert space with norm
\begin{equation}
  \|u\|_{H_w^s}^2 := \sum_{|\alpha| \le s} \int_{\R^m} |D_x^\alpha u(x)|^2 w(x)\dd x,
  \label{eq:smooth-dual-metric}
\end{equation}
and let $H_w^{-s}(\R^m) := (H_w^s(\R^m))^*$ be its continuous dual, endowed with the dual norm $\|\zeta\|_{-s,w} := \sup_{\|u\|_{H_w^s} \le 1} |\zeta(u)|$.

This choice of parameters plays two explicit roles:
\begin{enumerate}
  \item[\textup{(i)}] \emph{Sobolev embedding ($s > 4 + m/2$).} For every multi-index $|\alpha| \le 4$, the point-evaluation functional $E_x^\alpha u := D_x^\alpha u(x)$ is bounded on $H_w^s(\R^m)$ with $\|E_x^\alpha\|_{-s,w}^2 \le C w(x)^{-1}$, accommodating the fourth-order derivatives arising in the corrector remainder $\mathcal{R}_u$ below.
  \item[\textup{(ii)}] \emph{Integrability and dual moment budget ($p_* = m+3 > m$).} Because $w \in L^1(\R^m)$, mollifying bounded-Lipschitz test functions $\varphi$ ($\|\varphi\|_{\BL}\le 1$) at scale $\varrho \in (0,1]$ gives $\|\varphi - \varphi_\varrho\|_\infty \le C\varrho$ and $\|D^\alpha \varphi_\varrho\|_\infty \le C\varrho^{1-|\alpha|}$ ($1 \le |\alpha| \le s$), whence $\|\varphi_\varrho\|_{H_w^s} \le C\varrho^{1-s}$. For any two positive finite measures $\rho,\rho'$ with masses $\mathcal{Z},\mathcal{Z}'$, this yields the interpolation inequality
  \begin{equation}
    d_{\BL}(\rho,\rho') \le C\varrho(\mathcal{Z} + \mathcal{Z}') + C\varrho^{1-s}\|\rho - \rho'\|_{-s,w},
    \qquad 0 < \varrho \le 1.
    \label{eq:unnormalized-bl-interpolation}
  \end{equation}
  Furthermore, the inverse weight $w(X_t)^{-1} = (1+|X_t|^2)^{p_*/2}$ has finite $L^2$ expectation because $q \ge 2p_*$ in Assumption~\ref{ass:dissipative}\textup{(D4)}.
\end{enumerate}

\begin{lemma}
\label{lem:averaged-random-propagators}
For the localized averaged model and any $\cF_{j-1}^Z$-measurable $v \in H_w^{-s}(\R^m)$,
\begin{equation}
  \E^{\mathbb{Q}_\delta}\!\left[\|\bar Q_j^{R,*} v\|_{-s,w}^2 \;\middle|\; \cF_{j-1}^Z\right]
  \le e^{C_R\delta}\|v\|_{-s,w}^2.
  \label{eq:averaged-step-energy}
\end{equation}
Let $(e_\nu)_{\nu \ge 1}$ be an orthonormal basis of $H_w^s(\R^m)$. For a terminal block $k \le \lfloor T/\delta \rfloor$ and $1 \le j \le k$, define the backward random test functions
\[
  U_{j,\nu} := \bar Q_{j+1,Z_{j+1}}^R \cdots \bar Q_{k,Z_k}^R e_\nu,
  \qquad
  u_{j,\nu}(t) := \bar P_{\delta-t}^R\bigl(G_j^R(\cdot; Z_j) U_{j,\nu}\bigr),
  \qquad 0 \le t \le \delta,
\]
with the empty product interpreted as the identity. Then for every multi-index $|\alpha| \le 4$ and $x \in \R^m$,
\begin{equation}
  \sum_{\nu=1}^\infty \E_{Z_{j:k}}^{\mathbb{Q}_\delta} \bigl|D_x^\alpha u_{j,\nu}(t,x)\bigr|^2
  \le C_{T,R}\, w(x)^{-1}.
  \label{eq:random-test-trace}
\end{equation}
\end{lemma}

\begin{proof}
\emph{Step 1.} We first establish stability of the localized averaged semigroup on $H_w^s(\R^m)$ and $H_w^{-s}(\R^m)$.
By Lemma~\ref{lem:cutoff-experiment}, $\bar f^R, \bar a^R, \bar H^R \in C_b^{s+3}(\R^m)$ and $\bar a^R \succeq \underline{a}_R I_m$. Note that for the polynomial weight $w(x) = (1+|x|^2)^{-p_*/2}$, every ratio $|D^\beta w(x)|/w(x)$ ($|\beta|\ge 1$) is uniformly bounded on $\R^m$. Differentiating $\|u(t)\|_{H_w^s}^2$ along the Kolmogorov backward equation $\partial_t u = \bar\cL^R u$, integrating by parts in the highest-order diffusion terms, and absorbing lower-order commutators and weight derivatives via uniform ellipticity and Young's inequality yields
\[
  \frac{\dd}{\dd t}\|u(t)\|_{H_w^s}^2
  \le -c_R \sum_{|\alpha|=s} \|\nabla D^\alpha u(t)\|_{L^2(w)}^2 + C_R\|u(t)\|_{H_w^s}^2.
\]
Gronwall's inequality gives $\|\bar P_t^R\|_{H_w^s \to H_w^s} \le e^{C_R t}$, and by duality $\|\bar P_t^{R,*}\|_{H_w^{-s} \to H_w^{-s}} \le e^{C_R t}$.

\emph{Step 2.} We next bound the mean-square operator norm of the random Gaussian likelihood multiplier.
Define the whitened observation map $\tilde h_\delta(x) := R_\delta^{-1/2} \bar H^R(x) \in \R^d$. Because $\|R_\delta^{-1/2}\|_{\op} \le \sqrt{\delta/r_0}$ and $\bar H^R \in C_b^{s+3}(\R^m)$, multiplication by each component $\tilde h_{\delta,r}$ is a bounded linear operator on $H_w^s(\R^m)$ (and hence on $H_w^{-s}(\R^m)$) with operator norm at most $C_R\sqrt{\delta}$. On an auxiliary probability space carrying a $d$-dimensional standard Brownian motion $B$ on $[0,1]$, consider the linear $H_w^{-s}(\R^m)$-valued SDE
\[
  \dd v_t = \sum_{r=1}^d \tilde h_{\delta,r} v_t \dd B_t^r,
  \qquad v_0 = v.
\]
Its unique solution at $t=1$ is $v_1 = \exp\{\tilde h_\delta^\top B_1 - \frac{1}{2}|\tilde h_\delta|^2\} v_0$. By the Hilbert-space It\^o formula,
\[
  \E\|v_t\|_{-s,w}^2
  = \|v_0\|_{-s,w}^2 + \sum_{r=1}^d \int_0^t \E\|\tilde h_{\delta,r} v_u\|_{-s,w}^2 \dd u
  \le \|v_0\|_{-s,w}^2 + C_R \delta \int_0^t \E\|v_u\|_{-s,w}^2 \dd u,
\]
so Gronwall's inequality gives $\E\|v_1\|_{-s,w}^2 \le e^{C_R\delta}\|v_0\|_{-s,w}^2$. Under $\mathbb{Q}_\delta$, the vector $R_\delta^{-1/2} Z_j$ is $\mathcal{N}(0,I_d)$ independent of $\cF_{j-1}^Z$, so multiplication by $G_j^R(\cdot; Z_j)$ has the exact same law as $\exp\{\tilde h_\delta^\top B_1 - \frac{1}{2}|\tilde h_\delta|^2\}$. Combining this with Step~1 for $\bar P_\delta^{R,*}$ proves \eqref{eq:averaged-step-energy}.

\emph{Step 3.} We combine the preceding estimates to prove the trace bound \eqref{eq:random-test-trace}.
For any $|\alpha| \le 4$ and $x \in \R^m$, since $s > |\alpha| + m/2$, the local Sobolev embedding on $B_1(x)$ and comparability of $w$ on $B_1(x)$ yield
\[
  |D_x^\alpha u(x)|^2 \le C \|u\|_{H^s(B_1(x))}^2 \le C w(x)^{-1} \|u\|_{H_w^s}^2,
\]
so $\|E_x^\alpha\|_{-s,w}^2 \le C w(x)^{-1}$. By Parseval's identity in $H_w^s(\R^m)$,
\[
  \sum_{\nu=1}^\infty |D_x^\alpha u_{j,\nu}(t,x)|^2
  = \sum_{\nu=1}^\infty \bigl|\langle E_x^\alpha, u_{j,\nu}(t) \rangle\bigr|^2
  = \left\|\bar Q_k^{R,*} \cdots \bar Q_{j+1}^{R,*}\bigl(G_j^R(\cdot; Z_j)\bar P_{\delta-t}^{R,*} E_x^\alpha\bigr)\right\|_{-s,w}^2.
\]
Taking expectations iteratively over $Z_k, \dots, Z_j$ via \eqref{eq:averaged-step-energy} and using $k\delta \le T$ yields
\[
  \sum_{\nu=1}^\infty \E_{Z_{j:k}}^{\mathbb{Q}_\delta} |D_x^\alpha u_{j,\nu}(t,x)|^2
  \le e^{C_R T} \|E_x^\alpha\|_{-s,w}^2
  \le C_{T,R}\, w(x)^{-1},
\]
which establishes \eqref{eq:random-test-trace}.
\end{proof}

\begin{proposition}
\label{prop:forward-adapted-energy}
Under Assumption~\ref{ass:dissipative}, initialize both filters from $\bar\pi_0 = \ProjX \Pi_0^\varepsilon$. Then
\begin{equation}
  \max_{k\delta \le T}
  \left(\E^{\mathbb{Q}_\delta}\|\widehat\rho_k^{\varepsilon,R} - \bar\rho_k^R\|_{-s,w}^2\right)^{1/2}
  \le C_{T,R,M_0}\frac{\varepsilon}{\delta}.
  \label{eq:forward-unnormalized-energy-rate}
\end{equation}
Consequently,
\begin{equation}
  \max_{k\delta \le T} \E^{\mathbb{Q}_\delta}\, d_{\BL}(\widehat\rho_k^{\varepsilon,R}, \bar\rho_k^R)
  \le C_{T,R,M_0}\left(\frac{\varepsilon}{\delta}\right)^{1/s}.
  \label{eq:endpoint-signal-bl-rate}
\end{equation}
\end{proposition}

\begin{proof}
\emph{Step 1.} We first telescope the multiblock transition difference and expand each single-block error via Poisson correctors.
Fix $k \le \lfloor T/\delta \rfloor$ and abbreviate the localized signal by $(X_t, Y_t) := (X_t^{\varepsilon,R}, Y_t^{\varepsilon,R})$. The one-step joint backward propagator of the endpoint-observation model is $\widehat{\mathsf{Q}}_j F := P_\delta^{\varepsilon,R}(G_j^R F)$, where $G_j^R(x; Z_j)$ depends only on $x$. Define the cumulative unnormalized functional up to step $j-1$ by
\[
  \widehat\Lambda_{j-1}(F) := \Pi_0^\varepsilon\bigl(\widehat{\mathsf{Q}}_1 \cdots \widehat{\mathsf{Q}}_{j-1} F\bigr),
  \qquad \widehat\Lambda_0 := \Pi_0^\varepsilon.
\]
For any test function $\varphi \in H_w^s(\R^m)$, set $U_j := \bar Q_{j+1,Z_{j+1}}^R \cdots \bar Q_{k,Z_k}^R \varphi$. Using $\bar\pi_0 = \ProjX \Pi_0^\varepsilon$ and the telescoping identity for products of operators, we obtain
\begin{equation}
  (\widehat\rho_k^{\varepsilon,R} - \bar\rho_k^R)(\varphi)
  = \sum_{j=1}^k \widehat\Lambda_{j-1}\!\left(P_\delta^{\varepsilon,R}(G_j^R U_j) - \bar P_\delta^R(G_j^R U_j)\right),
  \label{eq:transition-telescoping}
\end{equation}
where $\bar P_\delta^R(G_j^R U_j)$ is lifted to $\R^{m+n}$ as a function independent of $y$.

Condition on the future observations $Z_{j:k}$ and set $u(t,x) := \bar P_{\delta-t}^R(G_j^R U_j)(x)$ for $t \in [0,\delta]$, which solves the backward Kolmogorov equation $\partial_t u + \bar\cL^R u = 0$ with terminal datum $u(\delta,x) = G_j^R(x; Z_j)U_j(x)$. Define the first-order Poisson corrector perturbation and its remainder by
\[
  K_u(t,x,y) := \chi_f^R(x,y) \cdot \nabla_x u(t,x) + \frac{1}{2}\chi_a^R(x,y) : D_x^2 u(t,x),
  \qquad
  \mathcal{R}_u := (\partial_t + \cL^{S,R}) K_u.
\]
Using $\cL^{\varepsilon,R} = \varepsilon^{-1}\cL_x^{F,R} + \cL^{S,R}$, $\partial_t u + \bar\cL^R u = 0$, and the Poisson equations \eqref{eq:signal-poisson}--\eqref{eq:diffusion-poisson}, the $O(1)$ fast-fluctuation terms cancel identically:
\[
  (\partial_t + \cL^{\varepsilon,R})(u + \varepsilon K_u)
  = (\cL^{S,R} - \bar\cL^R) u + \cL_x^{F,R} K_u + \varepsilon(\partial_t + \cL^{S,R}) K_u
  = \varepsilon \mathcal{R}_u.
\]
Applying It\^o's formula to $u(t,X_t) + \varepsilon K_u(t,X_t,Y_t)$ on $[0,\delta]$ (first for smooth test functions up to stopping times, then extending by density in $H_w^s(\R^m)$) gives the one-step transition expansion
\begin{align}
  P_\delta^{\varepsilon,R}(G_j^R U_j)(x,y) - \bar P_\delta^R(G_j^R U_j)(x)
  &= \varepsilon K_u(0,x,y)
   - \varepsilon \E_{x,y}\bigl[K_u(\delta, X_\delta, Y_\delta)\bigr] \nonumber\\
  &\quad + \varepsilon \E_{x,y}\!\left[\int_0^\delta \mathcal{R}_u(t, X_t, Y_t)\dd t\right].
  \label{eq:backward-corrector-transition}
\end{align}
Substituting $\partial_t u = -\bar\cL^R u$ into $\mathcal{R}_u$ yields derivatives $D_x^\alpha u$ only up to order $|\alpha| \le 4$, multiplied by corrector derivatives, bounded coefficients $(\bar f^R, \bar a^R, a^R)$, or at most one polynomially growing factor $f^R$. Applying Lemma~\ref{lem:cutoff-poisson} with $\eta = \Theta/32$ and absorbing polynomial growth into $V_{\Theta/32}$ gives, with $W = V_{\Theta/16}$,
\begin{equation}
  |K_u(t,x,y)| + |\mathcal{R}_u(t,x,y)|
  \le C_R W(y) \sum_{|\alpha| \le 4} |D_x^\alpha u(t,x)|.
  \label{eq:transition-remainder-budget}
\end{equation}

\emph{Step 2.} We next decouple past likelihoods from future random test evaluations to estimate each summand in \eqref{eq:transition-telescoping} in $L^2(\mathbb{Q}_\delta; H_w^{-s}(\R^m))$. Let $\widehat L_{j-1} := \prod_{i=1}^{j-1} G_i^R(X_{t_i}; Z_i)$ be the cumulative endpoint likelihood up to step $j-1$. Under $\mathbb{Q}_\delta$, conditional on the signal path, the observations $Z_1,\dots,Z_{j-1}$ are independent $\mathcal{N}(0,R_\delta)$, so direct Gaussian integration gives
\[
  \E_{Z_{1:j-1}}^{\mathbb{Q}_\delta}\bigl[\widehat L_{j-1}^{\,2}\bigr]
  = \exp\!\left\{\sum_{i=1}^{j-1} \bar H^R(X_{t_i})^\top R_\delta^{-1} \bar H^R(X_{t_i})\right\}
  \le \exp\!\left(\frac{M_R^2 T}{r_0}\right),
\]
where $M_R = \sup_x |\bar H^R(x)|$ from Lemma~\ref{lem:cutoff-experiment}. Moreover, under $\mathbb{Q}_\delta$, the future observations $Z_{j:k}$ are independent of both $Z_{1:j-1}$ and the signal path.

For either block endpoint $t \in \{t_{j-1}, t_j\}$ and basis function $e_\nu$ of $H_w^s(\R^m)$, denote the boundary term in \eqref{eq:backward-corrector-transition} by
\[
  B_{j,\nu} := \varepsilon\, \E_{\mathrm{sig}}\!\left[\widehat L_{j-1}\, K_{u_{j,\nu}}(t - t_{j-1}, X_t, Y_t)\right].
\]
For every finite $N \ge 1$, applying Jensen's inequality in $\E_{\mathrm{sig}}$, independence of $Z_{1:j-1}$ and $Z_{j:k}$, \eqref{eq:transition-remainder-budget}, and the trace bound \eqref{eq:random-test-trace} gives
\begin{align*}
  \sum_{\nu=1}^N \E^{\mathbb{Q}_\delta}|B_{j,\nu}|^2
  &\le \varepsilon^2\, \E_{\mathrm{sig}}\!\left[
    \E_{Z_{1:j-1}}^{\mathbb{Q}_\delta}\bigl[\widehat L_{j-1}^{\,2}\bigr]
    \sum_{\nu=1}^N \E_{Z_{j:k}}^{\mathbb{Q}_\delta}\bigl|K_{u_{j,\nu}}(t - t_{j-1}, X_t, Y_t)\bigr|^2
  \right] \\
  &\le C_{T,R}\,\varepsilon^2 \sup_{0 \le t \le T} \E_{\mathrm{sig}}\bigl[W(Y_t)^2 w(X_t)^{-1}\bigr].
\end{align*}
By Cauchy--Schwarz's inequality and Lemma~\ref{lem:localized-block-moments} (using $W^4 = V_{\Theta/4} \le V_\Theta$ and $q \ge 2p_*$),
\[
  \E\bigl[W(Y_t)^2 w(X_t)^{-1}\bigr]
  \le \bigl(\E[V_{\Theta/4}(Y_t)]\bigr)^{1/2}\bigl(\E[(1+|X_t|^2)^{p_*}]\bigr)^{1/2}
  \le C_{T,R,M_0}.
\]
The same argument with Cauchy--Schwarz's inequality in time bounds the $\mathcal{R}_u$-integral contribution by $C_{T,R,M_0}\varepsilon^2 \delta^2$.

\emph{Step 3.} We now sum over all blocks and pass from the dual Sobolev norm to the bounded-Lipschitz distance.
Letting $N \to \infty$ by monotone convergence and Parseval's identity shows that each $j$-th summand in \eqref{eq:transition-telescoping} has $L^2(\mathbb{Q}_\delta; H_w^{-s}(\R^m))$ norm at most $C_{T,R,M_0}\varepsilon$. Summing over $j=1,\dots,k \le T/\delta$ via Minkowski's inequality yields \eqref{eq:forward-unnormalized-energy-rate}. Finally, since $\E^{\mathbb{Q}_\delta}[\widehat\rho_k^{\varepsilon,R}(1)] = \E^{\mathbb{Q}_\delta}[\bar\rho_k^R(1)] = 1$, taking expectations in \eqref{eq:unnormalized-bl-interpolation} and setting $\varrho := (\varepsilon/\delta)^{1/s} \in (0,1]$ yields \eqref{eq:endpoint-signal-bl-rate}.
\end{proof}

\begin{proof}[Proof of Theorem~\ref{thm:filter-rate}]
\emph{Step 1.} We first verify well-posedness and the exit-probability bound for the original averaged signal.
Integrating \textup{(D3)} against $\mu_x$ and using the uniform invariant moment bound $\sup_x \mu_x(V_\eta) < \infty$ yields $\langle x, \bar f(x)\rangle \le -\kappa_s|x|^2 + C$ and $\operatorname{tr}\bar a(x) \le C$. By Lemma~\ref{lem:cutoff-poisson} and \textup{(D5)}, $\bar f$ and $\bar a^{1/2}$ are locally Lipschitz. Hence \eqref{eq:averaged-signal} has a unique nonexplosive strong solution $\bar X$, and the same It\^o--BDG--Gronwall argument as in Lemma~\ref{lem:signal-moments} gives $\E\sup_{t\le T}|\bar X_t|^q \le C_{T,M_0}$. For the exit times $\tau_R^\varepsilon := \inf\{t \ge 0 : |X_t^\varepsilon| \ge R\}$ and $\bar\tau_R := \inf\{t \ge 0 : |\bar X_t| \ge R\}$, Markov's inequality and Lemma~\ref{lem:signal-moments} yield
\begin{equation}
  \sup_{0 < \varepsilon \le 1} \Prob(\tau_R^\varepsilon \le T)
  + \bar\Prob(\bar\tau_R \le T)
  \le C_{T,M_0} R^{-q} =: \vartheta_T(R),
  \label{eq:exit-tail-bound}
\end{equation}
where $\bar\Prob$ denotes the law of the autonomous averaged signal initialized at $\bar\pi_0$.

\emph{Step 2.} We next bound the localization errors by pathwise coupling.
Fix $k \le \lfloor T/\delta \rfloor$ and let $Q_k := \mathbb{Q}_\delta|_{\cF_k^Z}$. Couple $(X^\varepsilon, Y^\varepsilon)$ with $(X^{\varepsilon,R}, Y^{\varepsilon,R})$ using the same initial state and driving Brownian motions $(W^X, W^Y)$, and couple their observations using the same Gaussian noise vectors $(\xi_j)$. Because $\psi_R(x) = x$ on $\{|x|\le R\}$, the original and localized trajectories and all block observables coincide identically on the event $\{\tau_R^\varepsilon > T\}$. Consequently, for any bounded measurable function $\Phi(x, z_{1:k})$ with $|\Phi| \le 1$,
\[
  \bigl|\E\bigl[\Phi(X_{t_k}^\varepsilon, Z_{1:k})\bigr]
       - \E\bigl[\Phi(X_{t_k}^{\varepsilon,R}, Z_{1:k}^R)\bigr]\bigr|
  \le 2\,\Prob(\tau_R^\varepsilon \le T).
\]
Writing the joint law of $(X_{t_k}^\varepsilon, Z_{1:k})$ as $\rho_k^\varepsilon(\dd x; z_{1:k})\,Q_k(\dd z_{1:k})$ and disintegrating with respect to $Q_k$ gives
\begin{equation}
  \E^{Q_k}\|\rho_k^\varepsilon - \rho_k^{\varepsilon,R}\|_{\TV}
  \le 2\,\Prob(\tau_R^\varepsilon \le T).
  \label{eq:full-posterior-coupling-bound}
\end{equation}
An identical coupling between $\bar X$ and $\bar X^R$ yields
\begin{equation}
  \E^{Q_k}\|\bar\rho_k - \bar\rho_k^R\|_{\TV}
  \le 2\,\bar\Prob(\bar\tau_R \le T).
  \label{eq:reduced-posterior-coupling-bound}
\end{equation}

\emph{Step 3.} We now assemble the preceding estimates to conclude.
Inserting the bounds \eqref{eq:full-posterior-coupling-bound}--\eqref{eq:reduced-posterior-coupling-bound}, Lemma~\ref{lem:block-endpoint-experiment}, and Proposition~\ref{prop:forward-adapted-energy} into the decomposition \eqref{eq:proof-decomposition-chain} yields
\[
  \E^{Q_k}\, d_{\BL}(\rho_k^\varepsilon, \bar\rho_k)
  \le C_{T,R,M_0}\!\left[r_{\varepsilon,\delta} + \left(\frac{\varepsilon}{\delta}\right)^{1/s}\right] + 2\vartheta_T(R).
\]
Applying Lemma~\ref{lem:bayes-joint-transfer} once to the unlocalized pair $(\rho_k^\varepsilon, \bar\rho_k)$ under $Q_k$ multiplies the right-hand side by $2$, establishing the first inequality in \eqref{eq:upgraded-main-theorem}. Simplifying $r_{\varepsilon,\delta} + (\varepsilon/\delta)^{1/s} \le C r_{\varepsilon,\delta}^{1/s}$ for $0 < \varepsilon \le \delta \le 1$ gives the second inequality. Taking $\varepsilon/\delta \to 0$ and $\delta \to 0$ at fixed $R$ proves \eqref{eq:iterated-limit-fixed-R}, and sending $R \to \infty$ completes the proof.
\end{proof}

\begin{remark}[Localized model and globally bounded coefficients]
\label{rem:localized-no-exit}
Combining Lemmas~\ref{lem:bayes-joint-transfer} and \ref{lem:block-endpoint-experiment} with Proposition~\ref{prop:forward-adapted-energy} directly yields the localized posterior bound
\begin{equation}
  \max_{k\delta \le T} \E^{\Prob^R}\, d_{\BL}(\pi_k^{\varepsilon,R}, \bar\pi_k^R)
  \le C_{T,R,M_0}\!\left[r_{\varepsilon,\delta} + \left(\frac{\varepsilon}{\delta}\right)^{1/s}\right]
  \le C_{T,R,M_0}\, r_{\varepsilon,\delta}^{1/s}.
  \label{eq:forward-full-filter-rate}
\end{equation}
In particular, if the original coefficients and observation map already satisfy global boundedness and uniform ellipticity in $x$, no radial truncation is needed and the exit tail $4\vartheta_T(R)$ in Theorem~\ref{thm:filter-rate} vanishes identically.
\end{remark}

\subsection{Continuous-Time Zakai Limit of the Reduced Filter}
\label{subsec:continuous-time-limit}

We now verify that as the block length $\delta \downarrow 0$, the discrete reduced filter driven by its own homogenized observations \eqref{eq:homogenized-obs-explicit} converges to the classical continuous-time nonlinear filter governed by the weak Zakai equation \cite{Zakai1969,BainCrisan2009}.

Let $\bar X$ solve \eqref{eq:averaged-signal} with initial law $\nu \in \cP(\R^m)$, and let $B$ be a $d$-dimensional standard Brownian motion independent of $\bar X$. Put $h := \bar H$, $t_k := k\delta$, and $C_\delta := \delta R_\delta = R_0 + \delta\widetilde{R}_\delta$. Couple the discrete and continuous observation channels on the same probability space via
\begin{align}
  \delta\bar Z_k^\delta
  &= \delta h(\bar X_{t_k}) + C_\delta^{1/2}(B_{t_k} - B_{t_{k-1}}),
  \label{eq:coupled-averaged-discrete-observation}\\
  \dd\bar{\mathcal{Y}}_t
  &= h(\bar X_t)\dd t + R_0^{1/2}\dd B_t,
  \qquad \bar{\mathcal{Y}}_0 = 0.
  \label{eq:continuous-homogenized-obs-channel}
\end{align}
Note that \eqref{eq:coupled-averaged-discrete-observation} reproduces the exact marginal law \eqref{eq:homogenized-obs-explicit} without requiring $R_0$ and $\widetilde{R}_\delta$ to commute. Define the discrete and continuous conditional distributions
\[
  \bar\pi_k^\delta := \operatorname{Law}(\bar X_{t_k} \mid \bar Z_1^\delta,\dots,\bar Z_k^\delta),
  \qquad
  \bar\pi_t^{\mathrm{ct}} := \operatorname{Law}\bigl(\bar X_t \mid \sigma(\bar{\mathcal{Y}}_u : 0 \le u \le t)\bigr),
\]
and let $\bar{\mathcal{Y}}^\delta$ denote the piecewise linear interpolation of the cumulative sums $\sum_{j=1}^k \delta\bar Z_j^\delta$ at the grid points $t_k$ (held constant after the final grid point up to $T$). The associated weak Zakai stochastic partial differential equation (SPDE) is
\begin{equation}
  \dd\bar\rho_t(\varphi)
  = \bar\rho_t(\bar\cL\varphi)\dd t
  + \bar\rho_t(\varphi\bar H)^\top R_0^{-1}\dd\bar{\mathcal{Y}}_t,
  \qquad \varphi \in C_c^\infty(\R^m),
  \qquad \bar\rho_0 = \nu.
  \label{eq:1d-zakai-spde-homogenized}
\end{equation}

\begin{theorem}[Continuous-time limit of the reduced filter]
\label{thm:discrete-continuous-filter}
Under Assumption~\ref{ass:dissipative} and \eqref{eq:observation-noise-scaling}, suppose $\nu(|x|^q) \le M_0$. For the coupled observation models \eqref{eq:coupled-averaged-discrete-observation}--\eqref{eq:continuous-homogenized-obs-channel}, every $T < \infty$, $R \ge R_{\mathrm{loc}}$, and $0 < \delta \le 1$,
\begin{equation}
  \max_{k\delta \le T} \E\|\bar\pi_k^\delta - \bar\pi_{t_k}^{\mathrm{ct}}\|_{\TV}
  \le C_{T,R,M_0}\sqrt{\delta} + C_{T,M_0} R^{-q}.
  \label{eq:discrete-continuous-filter-rate}
\end{equation}
In particular, $\max_{k\delta \le T} \E\|\bar\pi_k^\delta - \bar\pi_{t_k}^{\mathrm{ct}}\|_{\TV} \to 0$ (and hence also in expected $d_{\BL}$) as $\delta \downarrow 0$. Moreover, $\bar{\mathcal{Y}}^\delta \to \bar{\mathcal{Y}}$ uniformly in probability on $[0,T]$, and $\bar\pi_t^{\mathrm{ct}} = \bar\rho_t / \bar\rho_t(1)$ where $\bar\rho_t$ solves the weak Zakai equation \eqref{eq:1d-zakai-spde-homogenized}.
\end{theorem}

\begin{proof}
We first establish an $O(\sqrt{\delta})$ total-variation bound for the localized model (Steps~1--3), then remove the localization via pathwise coupling (Step~4), and finally identify the cumulative observation limit and weak Zakai SPDE (Step~5).

\emph{Step 1.} We first whiten the observations and express both localized filters via the Kallianpur--Striebel formula.
Fix $R \ge R_{\mathrm{loc}}$ and abbreviate $X := \bar X^R$ and $h := \bar H^R$ throughout Steps~1--3. By Lemma~\ref{lem:cutoff-experiment}, $X$ has bounded coefficients and $h$ is bounded and globally Lipschitz. Smoothness of the matrix square root on uniformly positive definite compact sets gives
\[
  \|C_\delta^{-1/2} - R_0^{-1/2}\|_{\op} + \|C_\delta^{1/2} - R_0^{1/2}\|_{\op} \le C\delta.
\]
Define the whitened maps $\tilde h := R_0^{-1/2} h$ and $\tilde h_\delta := C_\delta^{-1/2} h$. For $s \in (t_{j-1}, t_j]$, set
\[
  \beta_s := \tilde h(X_s),
  \qquad
  \beta_s^\delta := \tilde h_\delta(X_{t_j}),
  \qquad
  V_t := B_t + \int_0^t \beta_s\dd s,
  \qquad
  V_t^\delta := B_t + \int_0^t \beta_s^\delta\dd s,
\]
so that $V = R_0^{-1/2}\bar{\mathcal{Y}}^R$ and $\Delta V_j^\delta := V_{t_j}^\delta - V_{t_{j-1}}^\delta = C_\delta^{-1/2}\delta\bar Z_j^{\delta,R}$. Both $\beta$ and $\beta^\delta$ are uniformly bounded by a constant $M_R$. Using $\E|X_{t_j} - X_s|^2 \le C_R(t_j - s)$ on $(t_{j-1}, t_j]$ and $\|\tilde h_\delta - \tilde h\|_\infty \le C_R\delta$, we obtain
\begin{equation}
  \E\!\left[\int_0^{t_k} |\beta_s^\delta - \beta_s|^2\dd s\right]
  \le C_{T,R}\delta,
  \qquad k\delta \le T.
  \label{eq:continuous-likelihood-integrand-error}
\end{equation}

Let $X'$ be an independent copy of $X$ (also independent of $B$), with corresponding processes $\beta', \beta^{\delta\prime}$, and let $\E'$ denote expectation over $X'$ alone. At $t = t_k$, define the relative likelihoods
\begin{align*}
  L_t' &:= \exp\!\left\{\int_0^t (\beta_s')^\top \dd V_s - \frac{1}{2}\int_0^t |\beta_s'|^2\dd s\right\}, \\
  L_t^{\delta\prime} &:= \exp\!\left\{\sum_{j=1}^k \tilde h_\delta(X_{t_j}')^\top \Delta V_j^\delta
    - \frac{\delta}{2}\sum_{j=1}^k |\tilde h_\delta(X_{t_j}')|^2\right\}.
\end{align*}
Conditional on $(X, X')$, the integrands against $B$ are deterministic, so the right-endpoint evaluation in $L_t^{\delta\prime}$ involves no anticipating stochastic integral. By the Kallianpur--Striebel formula \cite[Theorem~3.22]{BainCrisan2009},
\[
  \bar\pi_t^{\mathrm{ct},R}(\varphi) = \frac{\E'[\varphi(X_t') L_t']}{\E'[L_t']},
  \qquad
  \bar\pi_k^{\delta,R}(\varphi) = \frac{\E'[\varphi(X_t') L_t^{\delta\prime}]}{\E'[L_t^{\delta\prime}]}.
\]

\emph{Step 2.} We next compare the discrete and continuous relative likelihoods in $L^2$ under the physical coupling.
For deterministic functions $\psi, \zeta \in L^2([0,t];\R^d)$, define the Gaussian exponential functional
\[
  F_t(\psi,\zeta) := \exp\!\left\{\int_0^t \psi_s^\top \dd B_s - \frac{1}{2}\|\psi\|_{L^2(0,t)}^2 + \langle \psi, \zeta \rangle_{L^2(0,t)}\right\}.
\]
We claim that on any $L^2(0,t)$-ball of radius $M_R\sqrt{T}$,
\begin{equation}
  \E_B\bigl|F_t(\psi,\zeta) - F_t(\tilde\psi,\tilde\zeta)\bigr|^2
  \le C_{T,R}\bigl(\|\psi - \tilde\psi\|_{L^2(0,t)}^2 + \|\zeta - \tilde\zeta\|_{L^2(0,t)}^2\bigr).
  \label{eq:gaussian-likelihood-lipschitz}
\end{equation}
Indeed, linearly interpolating $(\psi_\theta, \zeta_\theta) := \theta(\psi,\zeta) + (1-\theta)(\tilde\psi,\tilde\zeta)$ for $\theta \in [0,1]$ and differentiating in $\theta$ gives
\[
  \frac{\dd}{\dd\theta} F_t(\psi_\theta, \zeta_\theta)
  = F_t(\psi_\theta, \zeta_\theta)\!\left(\int_0^t \dot\psi^\top \dd B
    - \langle \psi_\theta, \dot\psi \rangle + \langle \dot\psi, \zeta_\theta \rangle + \langle \psi_\theta, \dot\zeta \rangle\right),
\]
where $\dot\psi = \psi - \tilde\psi$ and $\dot\zeta = \zeta - \tilde\zeta$. Applying Cauchy--Schwarz in $\E_B$, noting that $\E_B[F_t(\psi_\theta,\zeta_\theta)^4] \le C_{T,R}$ and $\E_B[(\int_0^t \dot\psi^\top \dd B)^4] = 3\|\dot\psi\|_{L^2(0,t)}^4$, and integrating over $\theta \in [0,1]$ proves \eqref{eq:gaussian-likelihood-lipschitz}.

Conditional on $(X,X')$, we have $L_t' = F_t(\beta', \beta)$ and $L_t^{\delta\prime} = F_t(\beta^{\delta\prime}, \beta^\delta)$, all with $L^2(0,t)$-norm bounded by $M_R\sqrt{T}$. Applying \eqref{eq:gaussian-likelihood-lipschitz} and \eqref{eq:continuous-likelihood-integrand-error} to both copies yields $\E\E'|L_t^{\delta\prime} - L_t'|^2 \le C_{T,R}\delta$. Setting $D_t^\delta := \E'|L_t^{\delta\prime} - L_t'|$, Jensen's inequality gives
\begin{equation}
  \E\bigl[(D_t^\delta)^2\bigr] \le C_{T,R}\delta.
  \label{eq:physical-likelihood-error}
\end{equation}

\emph{Step 3.} We now transfer the likelihood error to the normalized localized posteriors.
Let $\mathcal{Z}_t := \E'[L_t']$ and $\bar\beta_s := \E'[\beta_s']$ (which is deterministic with $|\bar\beta_s| \le M_R$). By Jensen's inequality for $\log$ and stochastic Fubini,
\[
  \log \mathcal{Z}_t \ge \E'[\log L_t']
  = \int_0^t \bar\beta_s^\top \dd B_s + \int_0^t \bar\beta_s^\top \beta_s \dd s - \frac{1}{2}\int_0^t \E'|\beta_s'|^2\dd s
  \ge \int_0^t \bar\beta_s^\top \dd B_s - \frac{3}{2} M_R^2 T.
\]
Hence $\E[\mathcal{Z}_t^{-2}] \le e^{3M_R^2 T} \E[\exp(-2\int_0^t \bar\beta_s^\top \dd B_s)] \le \exp(5M_R^2 T)$ uniformly on $[0,T]$. For any $\|\varphi\|_\infty \le 1$, the quotient algebra of Lemma~\ref{lem:bayes-joint-transfer} gives pointwise
\[
  \bigl|\bar\pi_k^{\delta,R}(\varphi) - \bar\pi_t^{\mathrm{ct},R}(\varphi)\bigr|
  \le \frac{2 D_t^\delta}{\mathcal{Z}_t}.
\]
Taking the supremum over $\|\varphi\|_\infty \le 1$, applying Cauchy--Schwarz, and using \eqref{eq:physical-likelihood-error} yields
\begin{equation}
  \max_{k\delta \le T} \E\|\bar\pi_k^{\delta,R} - \bar\pi_{t_k}^{\mathrm{ct},R}\|_{\TV}
  \le C_{T,R}\sqrt{\delta}.
  \label{eq:cutoff-discrete-continuous-filter-rate}
\end{equation}

\emph{Step 4.} We next remove the slow-state localization by pathwise coupling.
Couple $\bar X$ and $\bar X^R$ via identical initial conditions and driving Brownian motions, using the same independent observation Brownian motion $B$. By \eqref{eq:exit-tail-bound}, $p_R := \Prob(\sup_{t\le T}|\bar X_t| \ge R) \le C_{T,M_0} R^{-q}$. Outside this exit event, the original and localized signals and observations coincide identically on $[0,T]$.

In general, if two state--observation pairs $(S,O)$ and $(S^R,O^R)$ coincide outside an event of probability at most $p_R$, then their joint laws differ in total variation by at most $2p_R$. Disintegrating with respect to a common observation reference measure and applying the total-variation form of Lemma~\ref{lem:bayes-joint-transfer} gives $\E\|\pi(O) - \pi^R(O)\|_{\TV} \le 4p_R$. Since $\|\pi^R(O) - \pi^R(O^R)\|_{\TV} \le 2\,\mathbf{1}_{\{O \ne O^R\}}$, the triangle inequality gives $\E\|\pi(O) - \pi^R(O^R)\|_{\TV} \le 6p_R$. Applying this bound once to the discrete pair $(\bar\pi_k^\delta, \bar\pi_k^{\delta,R})$ and once to the continuous pair $(\bar\pi_{t_k}^{\mathrm{ct}}, \bar\pi_{t_k}^{\mathrm{ct},R})$, and combining with \eqref{eq:cutoff-discrete-continuous-filter-rate}, proves \eqref{eq:discrete-continuous-filter-rate} with exit remainder $12p_R \le C_{T,M_0} R^{-q}$.

\emph{Step 5.} Finally, we establish convergence of the cumulative observations and identify the weak Zakai SPDE.
At grid points $t_k$, \eqref{eq:coupled-averaged-discrete-observation} yields $\sum_{j=1}^k \delta\bar Z_j^\delta = \sum_{j=1}^k \delta h(\bar X_{t_j}) + C_\delta^{1/2} B_{t_k}$. Continuity of $t \mapsto h(\bar X_t)$ gives uniform convergence of the Riemann sums to $\int_0^\cdot h(\bar X_s)\dd s$, while $\|C_\delta^{1/2} - R_0^{1/2}\|_{\op} \le C\delta$ and Brownian continuity yield $\sup_{t\le T}|\bar{\mathcal{Y}}_t^\delta - \bar{\mathcal{Y}}_t| \to 0$ in probability.

For the weak Zakai equation, let $\mathbb{Q}_{\mathrm{ct}}$ be the reference law under which $\bar X$ retains its autonomous law and $\bar{\mathcal{Y}} = R_0^{1/2} B$ is independent of $\bar X$. Conditional on $\bar X$, the Girsanov density
\[
  L_t := \exp\!\left\{\int_0^t h(\bar X_s)^\top R_0^{-1}\dd\bar{\mathcal{Y}}_s
    - \frac{1}{2}\int_0^t h(\bar X_s)^\top R_0^{-1} h(\bar X_s)\dd s\right\}
\]
has conditional expectation one, and $|h(x)| \le C(1+|x|)$ implies $\E\int_0^T |h(\bar X_s)|^2\dd s < \infty$. For the localized model, applying It\^o's product rule to $\varphi(\bar X_t^R) L_t^R$ and taking $\mathbb{Q}_{\mathrm{ct}}$-conditional expectation given $\sigma(\bar{\mathcal{Y}}_u : u \le t)$ yields \eqref{eq:1d-zakai-spde-homogenized} with localized coefficients (the signal martingale projects to zero by independence of $\bar W$ and $B$).

To pass from the localized to the original model without requiring global second moments of $L_t$, note that pathwise coupling gives $\sup_{t\le T} \E^{\mathbb{Q}_{\mathrm{ct}}}\|\bar\rho_t^R - \bar\rho_t\|_{\TV} \le 2p_R$. For any fixed $\varphi \in C_c^\infty(\R^m)$ and sufficiently large $R$, $\bar\cL^R\varphi = \bar\cL\varphi$ and $\varphi\bar H^R = \varphi\bar H$ are bounded. Stopping when $\bar\rho_t^R(1) + \bar\rho_t(1)$ reaches a threshold $K$ (which occurs with probability at most $2/K$ since the sum is a nonnegative $\mathbb{Q}_{\mathrm{ct}}$-martingale starting at $2$), It\^o's isometry shows that the stochastic-integral difference vanishes as $R \to \infty$, while the drift integral vanishes in $L^1$. Letting $K \to \infty$ establishes \eqref{eq:1d-zakai-spde-homogenized}, and Bayes' formula gives $\bar\pi_t^{\mathrm{ct}} = \bar\rho_t / \bar\rho_t(1)$.
\end{proof}

\begin{remark}[Relation between Theorems~\ref{thm:filter-rate} and \ref{thm:discrete-continuous-filter}]
When $\widetilde{R}_\delta = 0$, the localized discrete propagator takes the Lie--Trotter splitting form
\[
  \bar Q_k^{R,*}\rho
  = \exp\!\left\{(\bar H^R)^\top R_0^{-1}\Delta\bar{\mathcal{Y}}_k^{\delta,R}
    - \frac{\delta}{2}(\bar H^R)^\top R_0^{-1}\bar H^R\right\}
    e^{\delta\bar\cL^{R,*}}\rho,
  \qquad \Delta\bar{\mathcal{Y}}_k^{\delta,R} := \delta\bar Z_k^{\delta,R}.
\]
Note that Theorem~\ref{thm:discrete-continuous-filter} drives the discrete reduced filter with observations generated by the \emph{reduced} model \eqref{eq:homogenized-obs-explicit}, whereas Theorem~\ref{thm:filter-rate} compares the full and reduced filters driven by the \emph{same} full-model observations \eqref{eq:realized-measure-observation}. Combining both results into a continuous-time limit for the full multiscale filter would additionally require comparing the observation laws of the full and reduced models.
\end{remark}

\section{Numerical Experiments}
\label{sec:examples}

This section evaluates the multiscale filter reduction in a controlled scalar nonlinear stochastic volatility model ($m=n=d=1$, $\ell_1 = c_1 = 1$) with Gaussian block observations. We examine (i) posterior convergence and computational speedup along the balanced scale path $\delta = \sqrt{\varepsilon}$, (ii) the U-shaped sampling-window tradeoff and the statistical role of the convexity adjustment, and (iii) slow-state point recovery and posterior interval calibration.

\subsection{Posterior Convergence and Computational Efficiency}
\label{sec:primary}

Consider the coupled slow--fast stochastic volatility signal
\begin{equation}
\begin{aligned}
  \dd X_t^\varepsilon &= \bigl[-X_t^\varepsilon - 0.25(X_t^\varepsilon)^3 + 4\tanh^2(Y_t^\varepsilon)\bigr]\dd t + 0.35\dd W_t^X, \\
  \dd Y_t^\varepsilon &= -\frac{8}{\varepsilon} Y_t^\varepsilon \dd t + \frac{0.55}{\sqrt{\varepsilon}}\dd W_t^Y,
\end{aligned}
\label{eq:signal}
\end{equation}
with independent standard Brownian motions $W^X, W^Y$, which satisfies all dissipativity, ellipticity, and smoothness conditions of Assumption~\ref{ass:dissipative}. Over each block $I_k = ((k-1)\delta, k\delta]$, the block observation is
\begin{equation}
  Z_k = \log\!\left(\frac{1}{\delta}\int_{I_k} e^{X_s^\varepsilon + Y_s^\varepsilon}\dd s\right)
        + \sqrt{\frac{R_0}{\delta}}\,\xi_k,
  \qquad \xi_k \overset{\mathrm{i.i.d.}}{\sim} \mathcal{N}(0,1).
  \label{eq:obs}
\end{equation}
The fast OU process has invariant law $\mu = \mathcal{N}(0, 0.55^2/16)$. Gauss--Hermite quadrature (40, 80, and 160 nodes, agreeing to within $2\times 10^{-17}$) yields
\[
  m_2 := \E_\mu[\tanh^2 Y] \simeq 0.018227,
  \qquad
  m_{\mathrm{c}} := \log \E_\mu[e^Y] = \frac{0.55^2}{32} = 0.009453125.
\]
The homogenized signal has drift $\bar f(x) = -x - 0.25 x^3 + 4m_2$, diffusion coefficient $0.35$, and effective observation map $\bar H(x) = x + m_{\mathrm{c}}$. The reduced filter alternates between solving the Fokker--Planck prediction PDE and applying the Gaussian Bayes update:
\begin{equation}
  \partial_t p = -\partial_x(\bar f p) + \frac{0.35^2}{2}\partial_{xx} p,
  \qquad
  \bar p_k(x) \propto \bar p_{k\mid k-1}(x)\exp\!\left(-\frac{(Z_k - x - m_{\mathrm{c}})^2}{2 R_0/\delta}\right).
  \label{eq:reduced}
\end{equation}

\paragraph{Numerical Discretization and Distance Evaluation.}
All runs use a target horizon $T = 2.5$ with $K = \lfloor T/\delta \rfloor$ observation blocks. Synthetic data are generated with $M_{\mathrm{data}} = \max\{64, \lceil 16\delta/\varepsilon \rceil\}$ substeps per block (using exact OU transitions for $Y^\varepsilon$, Euler--Maruyama updates for $X^\varepsilon$, and trapezoidal quadrature for the block integral). The full path-integral filter is approximated via a sequential Monte Carlo particle filter \cite{DelMoral2004} with $M_{\mathrm{PF}} = \max\{64, \lceil 8\delta/\varepsilon \rceil\}$ substeps per block, path-integral likelihood weighting, and systematic resampling whenever the effective sample size $\mathrm{ESS}_k := (\sum_{i=1}^N (\omega_k^i)^2)^{-1}$ drops below $N/2$. The reduced PDE \eqref{eq:reduced} is solved on a uniform grid of $2001$ nodes over $[-2.5, 2.5]$ with zero-flux boundary conditions, an exponentially fitted conservative finite-volume scheme, and implicit Euler stepping ($\lceil \delta/0.002 \rceil$ steps per block). At each scale, both filters share a common empirical prior generated from $1000$ burn-in endpoints (burn-in duration $20$, step size at most $\varepsilon/16$, initialized at $X=0$ and $Y\sim\mu$).

Posterior distances $d_{\BL}(\pi_k^\varepsilon, \bar\pi_k)$ are evaluated after each observation update and prior to resampling by solving the exact linear program on the ordered union $\{x_i\}$ of particle locations and PDE grid nodes (with discrete mass differences $a_i$ and primal/dual feasibility tolerances $10^{-9}$):
\[
  \max_{f, L} \sum_i a_i f_i
  \qquad\text{subject to}\qquad
  0 \le L \le 1,\quad
  |f_i| \le 1 - L,\quad
  |f_{i+1} - f_i| \le L(x_{i+1} - x_i).
\]

\paragraph{Balanced-Scale Convergence and Resolution Diagnostics.}
We fix $\delta = \sqrt{\varepsilon}$ and $R_0 = 0.004$ across ten scales $\varepsilon = 2^{-j}$ ($j = 5,\dots,14$). On each path, $d_{\BL}$ is averaged over eight equally spaced checkpoint block indices between $1$ and $K$ (rounded to integers without duplicates) to form the path-level metric $D_r$, from which we report the Monte Carlo mean $\widehat{D}$ and standard error (MC SE).

\begin{figure}[htbp]
  \centering
  \includegraphics[width=.9\textwidth]{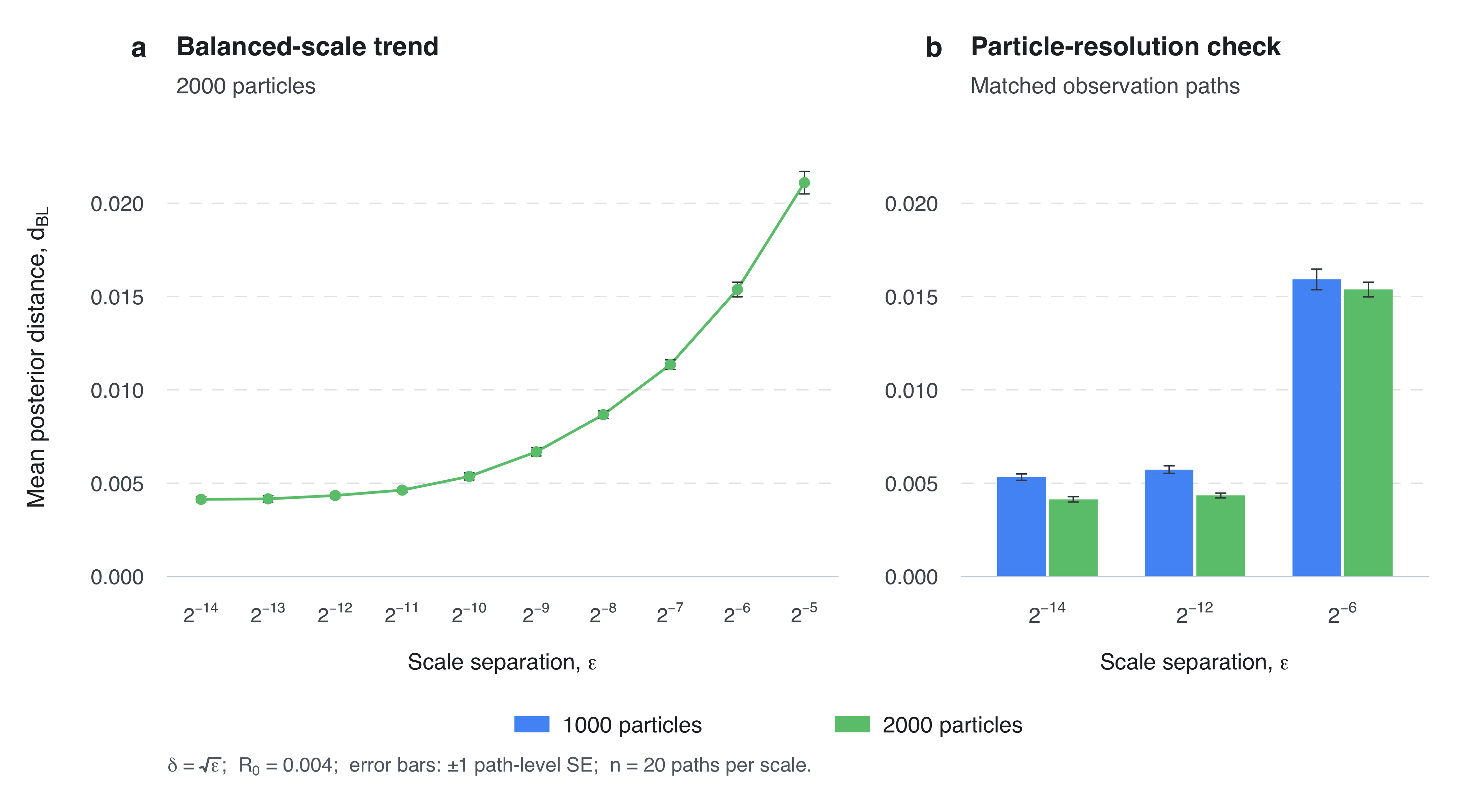}
  \caption{Posterior agreement and particle-resolution sensitivity
  along $\delta=\sqrt{\varepsilon}$, with $R_0=0.004$.
  \textbf{(a)} Mean bounded-Lipschitz distance $d_{\mathrm{BL}}$
  between the full filter's slow marginal and the reduced filter,
  using 2000 particles at ten scales
  $\varepsilon=2^{-j}$, $j=5,\ldots,14$.
  \textbf{(b)} Comparison of 1000- and 2000-particle full filters
  at $\varepsilon\in\{2^{-14},2^{-12},2^{-6}\}$,
  using shared observation paths and the same empirical prior
  within each scale.}
  \label{fig:primary}
\end{figure}

Table~\ref{tab:primary} and Figure~\ref{fig:primary} report the main $2000$-particle experiment ($20$ independent paths per scale, $1600$ checkpoint comparisons total). The mean bounded-Lipschitz distance decreases monotonically from $0.021100$ at $\varepsilon = 2^{-5}$ to $0.004133$ at $\varepsilon = 2^{-14}$ (an $80.4\%$ reduction).

\begin{table}[!htbp]
\centering
\caption{Balanced-scale posterior convergence ($\delta = \sqrt{\varepsilon}$, $R_0 = 0.004$) using $2000$ particles across $20$ independent paths per scale.}
\label{tab:primary}
\begin{tabular}{rrrrr}
\toprule
$\varepsilon$ & $\delta$ & $K$ & Mean $d_{\BL}$ & MC SE \\
\midrule
$2^{-5}$  & 0.176777 & 14  & 0.021100 & 0.000602 \\
$2^{-6}$  & 0.125000 & 20  & 0.015376 & 0.000393 \\
$2^{-7}$  & 0.088388 & 28  & 0.011353 & 0.000254 \\
$2^{-8}$  & 0.062500 & 40  & 0.008669 & 0.000208 \\
$2^{-9}$  & 0.044194 & 56  & 0.006678 & 0.000220 \\
$2^{-10}$ & 0.031250 & 80  & 0.005361 & 0.000189 \\
$2^{-11}$ & 0.022097 & 113 & 0.004627 & 0.000111 \\
$2^{-12}$ & 0.015625 & 160 & 0.004339 & 0.000130 \\
$2^{-13}$ & 0.011049 & 226 & 0.004159 & 0.000173 \\
$2^{-14}$ & 0.007812 & 320 & 0.004133 & 0.000142 \\
\bottomrule
\end{tabular}
\end{table}

To diagnose the fine-scale flattening in Figure~\ref{fig:primary}(a), we compare against a baseline $1000$-particle sweep ($50$ paths per scale, $4000$ checkpoints total), in which the mean distance falls from $0.021073$ at $\varepsilon = 2^{-5}$ to $0.005395$ at $\varepsilon = 2^{-14}$ ($74.4\%$ reduction). Across those $1000$-particle runs, mean pre-resampling ESS rises from $544.4$ at $\varepsilon = 2^{-5}$ to $704.1$ at $2^{-10}$, $736.4$ at $2^{-12}$, and $753.6$ at $2^{-14}$, while the resampling frequency drops from $41.4\%$ to $11.9\%$, $6.4\%$, and $3.5\%$. Because $R_\delta = R_0/\delta$ increases as $\delta$ shrinks, individual block updates become gentler and weight degeneracy does not occur; nevertheless, finite-particle propagation error accumulates over the growing number of blocks ($K$ rises from $14$ to $320$). One-at-a-time numerical sensitivity checks on $20$ paired paths per scale confirm that doubling the particle count from $1000$ to $2000$ reduces the mean distance by $0.001389$ at $\varepsilon = 2^{-12}$ and by $0.001193$ at $\varepsilon = 2^{-14}$ (both exceeding their paired standard errors). By contrast, refining the reduced PDE grid and time step alters the mean by at most $0.000027$, while refining the particle-filter and data substeps at $\varepsilon = 2^{-14}$ shifts the mean by $0.000280$ and $0.000331$, respectively. Finite-particle Monte Carlo error therefore accounts for much of the fine-scale plateau.

Table~\ref{tab:runtime} illustrates the computational advantage of the reduced filter. Because the full particle filter must resolve the fast time scale $\varepsilon$, its runtime per path grows from $0.0436\text{ s}$ at $\varepsilon = 2^{-5}$ to $16.3779\text{ s}$ at $\varepsilon = 2^{-14}$, whereas the reduced PDE filter remains nearly $\varepsilon$-independent ($0.0203$--$0.0296\text{ s}$), yielding a speedup exceeding $550\times$ at the finest scale.

\begin{table}[!htbp]
\centering\small
\caption{Per-path runtime comparison between the full particle filter and the reduced PDE filter along $\delta = \sqrt{\varepsilon}$.}
\label{tab:runtime}
\begin{tabular}{rrrrrr}
\toprule
$\varepsilon$ & Blocks & PF substeps/block & Full PF (s) & Reduced PDE (s) & Speedup ratio \\
\midrule
$2^{-5}$  & 14  & 64   & 0.0436  & 0.0203 & 2.2 \\
$2^{-10}$ & 80  & 256  & 0.9925  & 0.0203 & 46.7 \\
$2^{-14}$ & 320 & 1024 & 16.3779 & 0.0296 & 552.9 \\
\bottomrule
\end{tabular}
\end{table}

\subsection{Sampling-Window Tradeoff and the Role of the Convexity Adjustment}
\label{sec:unbalanced}

To examine the U-shaped dependence of posterior error on the observation window $\delta$ and isolate the systematic shift removed by the convexity adjustment, we retain the slow dynamics in \eqref{eq:signal} and adopt a fast OU component with larger stationary variance:
\begin{equation}
  \dd Y_t^\varepsilon = -\frac{\kappa_f}{\varepsilon} Y_t^\varepsilon \dd t
                        + \frac{\nu_f}{\sqrt{\varepsilon}}\dd W_t^Y,
  \qquad \varepsilon = 2^{-9},\quad \kappa_f = 1,\quad \nu_f = 0.8.
  \label{eq:tradeoff-fast-model}
\end{equation}
Here the stationary fast variance is $v_f := \nu_f^2/(2\kappa_f) = 0.32$, the fast correlation time is $\tau_f = \varepsilon/\kappa_f$, and the theoretical convexity shift is $m_{\mathrm{c}} = v_f/2 = 0.16$. The averaged slow drift is $-x - 0.25 x^3 + 4m_2$ with $m_2 = \E_{\mathcal{N}(0,0.32)}[\tanh^2 Y] \simeq 0.206138$. Fixing $R_0 = 0.004$ and identical averaged dynamics, we compare the convexity-adjusted observation map $H_{\mathrm{corr}}(x) := x + 0.16$ against the uncorrected map $H_{\mathrm{uncorr}}(x) := x$ across seven dyadic windows
\begin{equation}
  \delta \in \{2^{-7}, 2^{-6}, 2^{-5}, 2^{-4}, 2^{-3}, 2^{-2}, 2^{-1}\}.
  \label{eq:unbalanced}
\end{equation}
The geometric-mean scale $\sqrt{\varepsilon} \simeq 0.0442$ lies between $2^{-5}$ and $2^{-4}$.

We simulate $50$ independent signal paths on $[0,2]$ starting from $X_0 = 0.75$ and $Y_0 \sim \mathcal{N}(0,0.32)$, coupling observation errors across windows via $\sqrt{R_0}(B_{t_k} - B_{t_{k-1}})/\delta$ (variance $R_0/\delta$), where $B$ is a standard Brownian motion independent of the signal. For each path $r$, we evaluate the four-checkpoint mean distances
\[
  D_r^{\mathrm{corr}}(\delta) := \frac{1}{4}\sum_{i=1}^4 d_{\BL}\bigl(\pi_{r,t_i}^\varepsilon, \bar\pi_{r,t_i}^{\mathrm{corr}}\bigr),
  \qquad
  D_r^{\mathrm{uncorr}}(\delta) := \frac{1}{4}\sum_{i=1}^4 d_{\BL}\bigl(\pi_{r,t_i}^\varepsilon, \bar\pi_{r,t_i}^{\mathrm{uncorr}}\bigr),
\]
at $t_i \in \{0.5, 1.0, 1.5, 2.0\}$. The full particle filter uses $4000$ particles (exact OU transitions, slow step at most $\varepsilon/8$, trapezoidal quadrature, systematic resampling at $\mathrm{ESS} < N/2$), signal generation uses step $\varepsilon/32$, and the reduced PDE uses $1001$ nodes on $[-1.5, 3.0]$ with implicit Euler step at most $0.001$ (LP tolerances $10^{-8}$).

\begin{figure}[!htbp]\centering
\includegraphics[width=.65\textwidth]{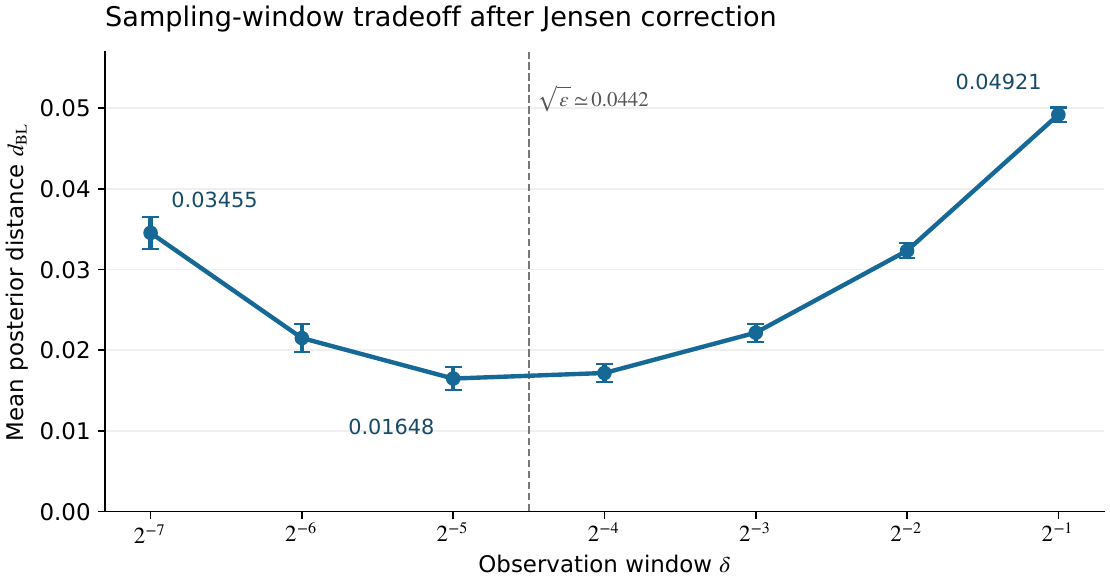}
\caption{Sampling-window tradeoff after Jensen correction: the full-versus-reduced posterior distance first decreases and then increases with $\delta$. The uncorrected comparison is reported in Table~\ref{tab:unbalanced}.}\label{fig:unbalanced}
\end{figure}

\paragraph{The U-Shaped Posterior-Distance Curve.}
As shown in Table~\ref{tab:unbalanced} and Figure~\ref{fig:unbalanced}, the convexity-adjusted posterior distance decreases from $0.034553$ at $\delta = 2^{-7}$ to a minimum of $0.016484$ at $\delta = 2^{-5}$ (near $\sqrt{\varepsilon}$), and then rises to $0.049206$ at $\delta = 2^{-1}$. At $\delta = 2^{-5}$, incorporating the convexity adjustment reduces the posterior discrepancy from $0.078210$ to $0.016484$ (a $79\%$ reduction). Relative to $\delta = 2^{-5}$, paired $95\%$ confidence intervals for the short- and long-window increases are $[0.016654, 0.019484]$ and $[0.031027, 0.034418]$, while the difference between $2^{-4}$ and $2^{-5}$ has $95\%$ interval $[-0.000046, 0.001429]$, confirming a well-resolved interior valley near $\sqrt{\varepsilon}$.

\begin{table}[!htbp]
\centering\small
\caption{Sampling-window tradeoff at $\varepsilon = 2^{-9}$, $\kappa_f = 1$, $\nu_f = 0.8$, and $R_0 = 0.004$ ($T=2$, $4000$ particles, $50$ paired paths; MC SE in parentheses).}
\label{tab:unbalanced}
\begin{tabular}{rrrr}
\toprule
$\delta$ & Blocks & Convexity-adjusted & Without correction \\
\midrule
0.007812 & 256 & $0.034553\;(0.000990)$ & $0.053663\;(0.000774)$ \\
0.015625 & 128 & $0.021507\;(0.000887)$ & $0.067944\;(0.000794)$ \\
0.031250 & 64  & $0.016484\;(0.000709)$ & $0.078210\;(0.000892)$ \\
0.062500 & 32  & $0.017175\;(0.000557)$ & $0.085708\;(0.000986)$ \\
0.125000 & 16  & $0.022175\;(0.000552)$ & $0.093627\;(0.001280)$ \\
0.250000 & 8   & $0.032342\;(0.000468)$ & $0.102521\;(0.001666)$ \\
0.500000 & 4   & $0.049206\;(0.000443)$ & $0.113824\;(0.001931)$ \\
\bottomrule
\end{tabular}
\end{table}

\begin{figure}[!htbp]\centering
\includegraphics[width=.7\textwidth]{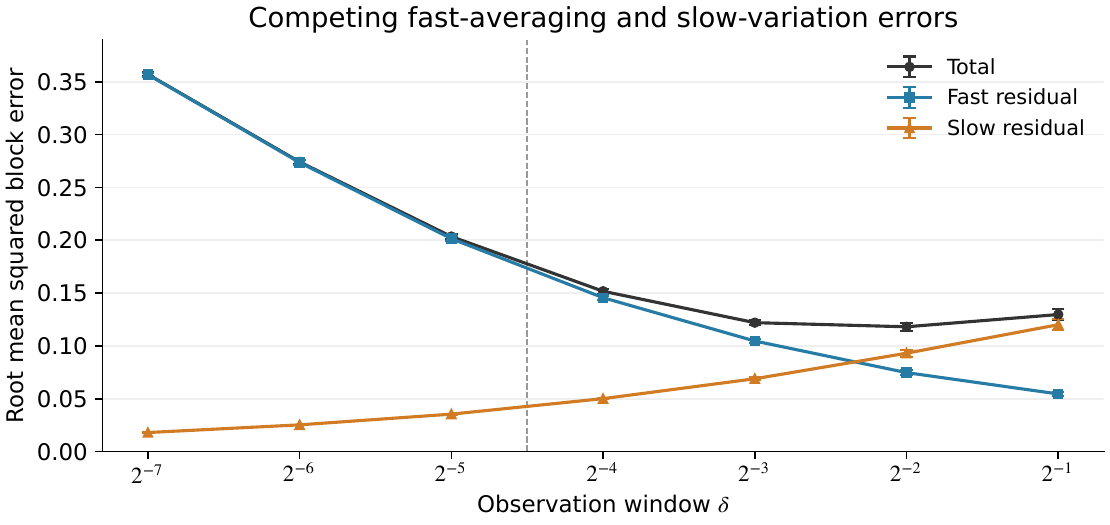}
\caption{Competing block errors underlying the window tradeoff in Figure~\ref{fig:unbalanced}. Root mean squared corrected block error and its fast and slow residuals are evaluated on 256 independent paths. The total includes the signed cross term. }\label{fig:tradeoff-support}
\end{figure}

\paragraph{Competing Fast and Slow Block Errors.}
To decompose the two arms of the U-shaped curve, we evaluate on $256$ independent paths the exact block-residual identity
\begin{align*}
  L_k^X &:= \log\!\left(\frac{1}{\delta}\int_{I_k} e^{X_s^\varepsilon}\dd s\right), \\
  \eta_k^{\mathrm{f}} &:= H_k^\varepsilon - L_k^X - m_{\mathrm{c}},
  \qquad
  \eta_k^{\mathrm{s}} := L_k^X - X_{t_k}^\varepsilon,
  \qquad
  H_k^\varepsilon - (X_{t_k}^\varepsilon + m_{\mathrm{c}}) = \eta_k^{\mathrm{f}} + \eta_k^{\mathrm{s}}.
\end{align*}
Here $\eta_k^{\mathrm{f}}$ measures incomplete fast ergodic averaging (weighted by the slow path), while $\eta_k^{\mathrm{s}}$ measures slow-trajectory variation over $I_k$. Figure~\ref{fig:tradeoff-support} confirms that the root mean square (RMS) of $\eta_k^{\mathrm{f}}$ decreases with $\delta$ while that of $\eta_k^{\mathrm{s}}$ increases with $\delta$, producing the interior minimum.

\begin{figure}[!htbp]\centering
\includegraphics[width=.8\textwidth]{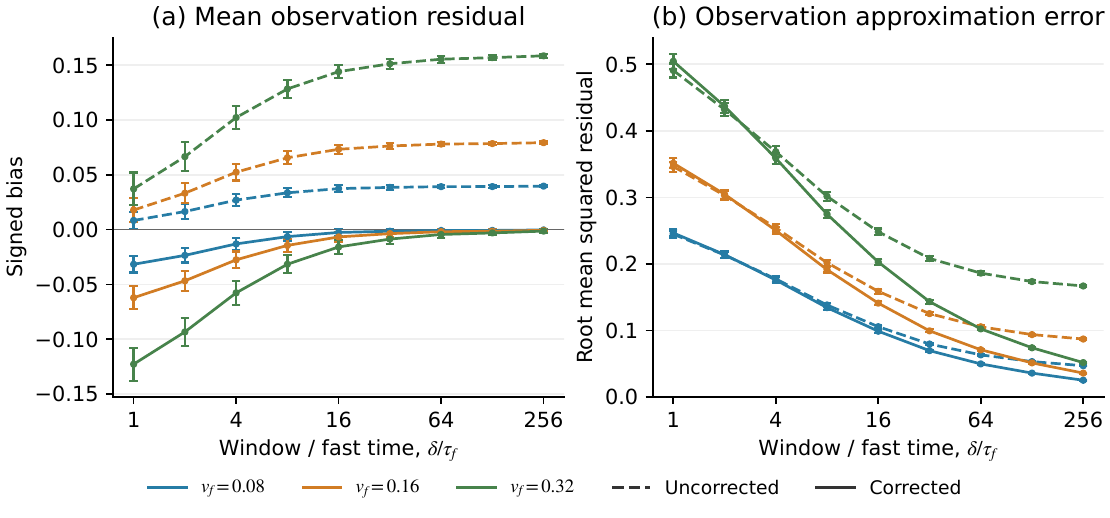}
\caption{Frozen-slow-state control with 4096 independent stationary OU paths. Dashed curves use residual $Q_\delta$; solid curves use $Q_\delta-v_f/2$. Colors distinguish invariant fast variances. (a) Signed mean residual. (b) Root mean squared residual. Error bars are pointwise 95\% path-level intervals; mean-square intervals are transformed by the square root. Paths are shared across variances and nested windows.}\label{fig:frozen-jensen}
\end{figure}

\paragraph{Frozen-Slow-State Control for the Convexity Shift.}
Freezing the slow state at $x$ isolates the block log-variance residual $Q_\delta := \log(\delta^{-1}\int_0^\delta e^{Y_s}\dd s)$ for a stationary fast OU process. Although $\E[Y_s] = 0$, ergodicity gives $Q_\delta \to v_f/2$ as $\delta/\tau_f \to \infty$. Simulating $4096$ stationary OU paths across $v_f \in \{0.08, 0.16, 0.32\}$ and $\delta/\tau_f \in \{1,2,4,\dots,256\}$ (Figure~\ref{fig:frozen-jensen}), we find that for $v_f = 0.32$ at $\delta/\tau_f = 256$, subtracting $v_f/2$ reduces the mean residual from $0.15860$ to $-0.00140$ and the RMS error from $0.16677$ to $0.05159$. At $\delta/\tau_f = 1$, by contrast, fast averaging has not yet occurred (corrected mean $-0.12287$, RMS $0.50456$ vs.\ $0.49078$ uncorrected), illustrating why the asymptotic shift requires $\delta \gg \tau_f$.

\paragraph{Consequences for Slow-State Inference.}
Omitting the positive convexity shift misattributes fast volatility variance to the slow factor, inducing upward location bias. Table~\ref{tab:jensen-state} summarizes slow-state posterior-mean bias, RMSE, equal-tailed $90\%$ credible-interval coverage, and interval width across the $50$ paths. At $\delta = 2^{-5}$, the convexity adjustment eliminates the signed bias ($0.1063 \to -0.0112$), lowers RMSE from $0.1636$ to $0.1324$ (paired improvement $0.0312$, $95\%$ interval $[0.0162, 0.0462]$), and raises $90\%$ coverage from $0.750$ to $0.865$ (paired gain $0.115$, $95\%$ interval $[0.0489, 0.1811]$) with negligible change in interval width ($0.4114$ vs.\ $0.4165$). A supplementary scaling ablation $H_\alpha(x) = x + \alpha m_{\mathrm{c}}$ at $\delta = 2^{-5}$ yields mean biases $(0.1063, 0.0480, -0.0112, -0.0714)$ and RMSEs $(0.1636, 0.1376, 0.1324, 0.1511)$ for $\alpha \in \{0, 1/2, 1, 3/2\}$, confirming that the theoretical value $\alpha = 1$ achieves the lowest error among the tested shifts.

\begin{table}[!htbp]
\centering\footnotesize
\caption{Slow-state inference in the strong-correction model ($50$ paired paths, four checkpoints per path; MC SE in parentheses, nominal coverage $0.90$).}
\label{tab:jensen-state}
\setlength{\tabcolsep}{3.5pt}
\begin{tabular}{rlrrrr}
\toprule
$\delta$ & Map & Signed bias & RMSE & Coverage & Width \\
\midrule
$2^{-7}$ & Uncorrected & $0.0738\;(0.0110)$  & $0.1473\;(0.0085)$ & $0.8300\;(0.0298)$ & $0.4235\;(0.0007)$ \\
$2^{-7}$ & Corrected   & $-0.0427\;(0.0109)$ & $0.1399\;(0.0066)$ & $0.8500\;(0.0202)$ & $0.4286\;(0.0006)$ \\
$2^{-5}$ & Uncorrected & $0.1063\;(0.0110)$  & $0.1636\;(0.0092)$ & $0.7500\;(0.0343)$ & $0.4114\;(0.0007)$ \\
$2^{-5}$ & Corrected   & $-0.0112\;(0.0108)$ & $0.1324\;(0.0069)$ & $0.8650\;(0.0228)$ & $0.4165\;(0.0006)$ \\
$2^{-1}$ & Uncorrected & $0.1412\;(0.0104)$  & $0.1977\;(0.0090)$ & $0.4750\;(0.0329)$ & $0.2626\;(0.0003)$ \\
$2^{-1}$ & Corrected   & $0.0035\;(0.0103)$  & $0.1446\;(0.0068)$ & $0.5750\;(0.0359)$ & $0.2644\;(0.0002)$ \\
\bottomrule
\end{tabular}
\end{table}

\paragraph{Numerical Sensitivity Checks.}
Resolution checks at $\delta \in \{2^{-7}, 2^{-5}, 2^{-1}\}$ verify stability of these conclusions: doubling the particle count to $8000$ ($20$ paired paths) shifts mean distances by $-0.000986$, $-0.000495$, and $0.000132$, halving the particle integration step ($10$ paired paths) shifts them by $-0.001144$, $0.000109$, and $0.000034$, and refining the reduced PDE to $2001$ nodes and step $0.0005$ changes RMSE by less than $7\times 10^{-6}$ and coverage by at most $0.005$.

\subsection{Slow-State Recovery and Posterior Interval Calibration}
\label{sec:recovery-diagnostics}

Finally, we evaluate latent slow-state recovery along the balanced window path $\delta = \sqrt{\varepsilon}$ at $\varepsilon \in \{2^{-6}, 2^{-9}, 2^{-12}\}$ in the baseline model \eqref{eq:signal}--\eqref{eq:obs} with $R_0 = 0.004$, using $50$ independent paths per scale (distinct from the primary convergence runs). At every observation time, both filters compute posterior means and equal-tailed $90\%$ credible intervals (weighted particle quantiles for the full filter; grid-mass quantiles for the reduced filter).

\begin{figure}[htbp]\centering
\includegraphics[width=0.7\textwidth]{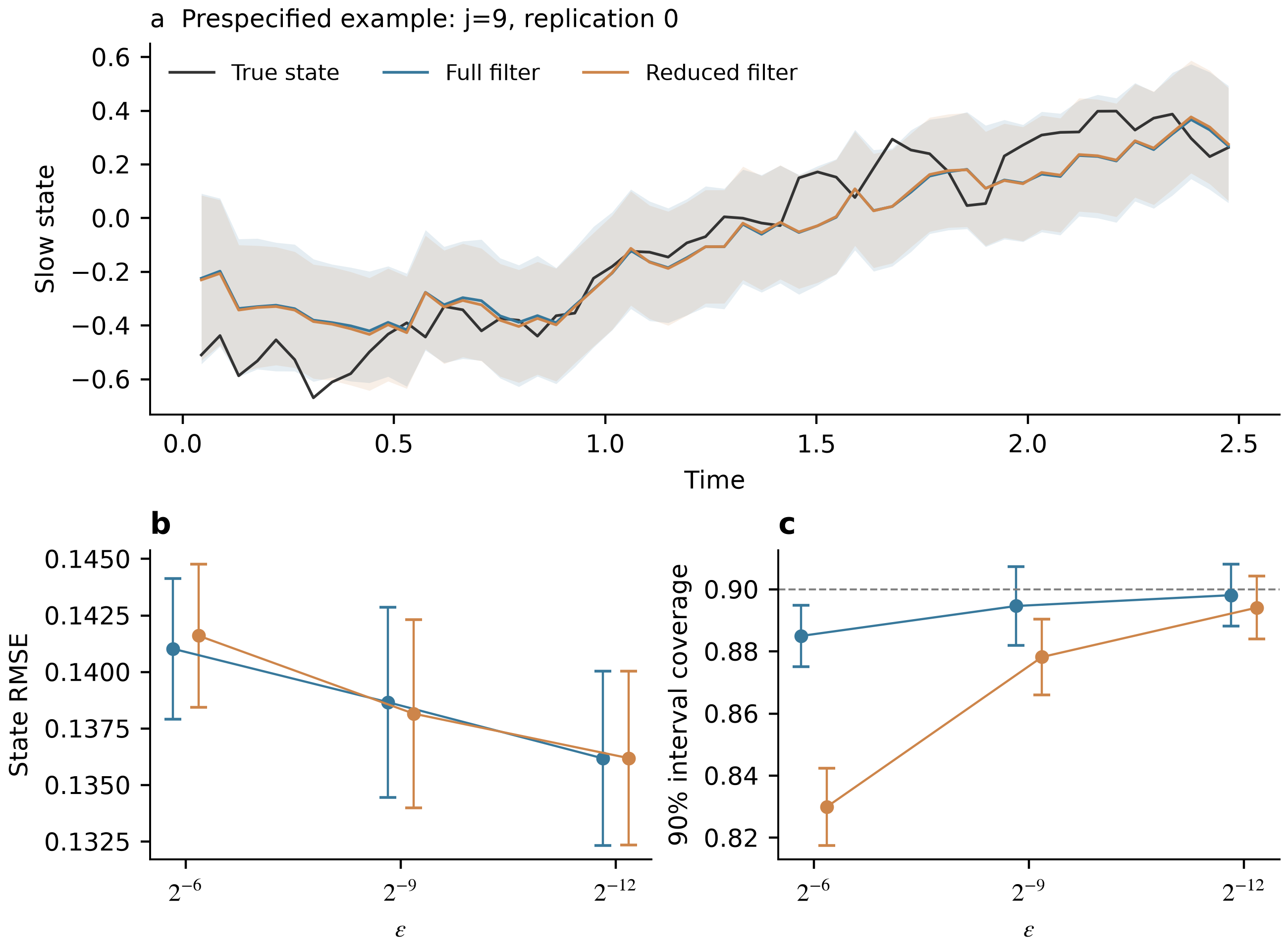}
\caption{State-recovery example and aggregate metrics. (a) The prespecified path with $j=9$, $a=1/2$, and replication 0. Lines show the true state and posterior means; shaded bands show the respective equal-tailed 90\% posterior intervals. This example does not represent aggregate coverage. (b--c) RMSE and coverage from Panel B of Table~\ref{tab:filter-diagnostics}, with $\pm1$ path-level MC SE and $n=50$ paths per scale. The dashed coverage line marks the nominal 90\% level.}
\label{fig:recovery}
\end{figure}

Table~\ref{tab:filter-diagnostics} and Figure~\ref{fig:recovery} summarize the results. At the coarse scale $\varepsilon = 2^{-6}$, the two filters achieve nearly identical RMSE ($0.1410$ vs.\ $0.1416$), yet the reduced filter's $90\%$ interval coverage is noticeably lower ($0.830$ vs.\ $0.885$) due to slightly narrower intervals ($0.3913$ vs.\ $0.4547$). As $\varepsilon$ decreases, both point accuracy and uncertainty calibration converge: at $\varepsilon = 2^{-9}$, RMSEs are $0.1387$ and $0.1382$ (coverages $0.895$ and $0.878$), and at $\varepsilon = 2^{-12}$, both RMSEs round to $0.1362$ with coverages $0.898$ and $0.894$. Thus, posterior interval calibration provides a stricter diagnostic of filter fidelity than point RMSE alone, and both align closely once the scale separation is sufficiently strong.

\begin{table}[!htbp]
\centering\small
\caption{Slow-state recovery along $\delta = \sqrt{\varepsilon}$ with $R_0 = 0.004$ ($50$ independent paths per scale; MC SE in parentheses, nominal coverage $0.90$).}
\label{tab:filter-diagnostics}
\begin{tabular}{rlrrrr}
\toprule
$\varepsilon$ & Filter & RMSE & Signed bias & 90\% coverage & Interval width \\
\midrule
$2^{-6}$  & Full    & $0.1410\;(0.0031)$ & $-0.0091\;(0.0079)$ & $0.8850\;(0.0099)$ & $0.4547\;(0.0006)$ \\
$2^{-6}$  & Reduced & $0.1416\;(0.0032)$ & $-0.0091\;(0.0075)$ & $0.8300\;(0.0125)$ & $0.3913\;(0.0001)$ \\
$2^{-9}$  & Full    & $0.1387\;(0.0042)$ & $0.0008\;(0.0077)$  & $0.8946\;(0.0127)$ & $0.4528\;(0.0006)$ \\
$2^{-9}$  & Reduced & $0.1382\;(0.0042)$ & $0.0004\;(0.0074)$  & $0.8782\;(0.0123)$ & $0.4313\;(0.0002)$ \\
$2^{-12}$ & Full    & $0.1362\;(0.0039)$ & $0.0000\;(0.0082)$  & $0.8981\;(0.0099)$ & $0.4543\;(0.0006)$ \\
$2^{-12}$ & Reduced & $0.1362\;(0.0038)$ & $0.0002\;(0.0082)$  & $0.8941\;(0.0101)$ & $0.4469\;(0.0002)$ \\
\bottomrule
\end{tabular}
\end{table}

\section{Discussion and Concluding Remarks}
\label{sec:conclusion}

We have developed a quantitative dimension-reduction framework for nonlinear filtering of multiscale stochastic volatility diffusions observed through logarithmic block realized measures. Two practical implications stand out for high-frequency volatility inference. First, because realized variance aggregates spot variance prior to taking logarithms, the reduced observation map must incorporate the state-dependent convexity adjustment $\log\mu_x(e^{c_r^\top Y})$; omitting this term misattributes fast-scale variance to the persistent volatility factor and degrades both point estimation and credible-interval calibration. Second, balancing fast ergodic averaging against slow-state displacement selects an intermediate sampling window $\Delta^* \asymp \sqrt{\tau_f\tau_s}$, at which the reduced PDE filter matches the slow marginal of the full path-integral filter while achieving more than two orders of magnitude computational speedup at fine scales.

Several concrete theoretical questions remain open for future investigation. First, quantifying the explicit dependence of the localized constant $C_{T,R,M_0}$ on the truncation radius $R$ would yield a global convergence rate without a separate exit remainder. Second, combining our single-block error estimates with nonlinear filter stability could extend the approximation from finite horizons to time-uniform bounds as $T \to \infty$. Third, establishing a direct comparison between the full and reduced observation laws would unify Theorems~\ref{thm:filter-rate} and \ref{thm:discrete-continuous-filter} into a single continuous-time Zakai limit for the original multiscale filter, while incorporating correlated slow--fast noises and path-dependent quarticity covariances (as in Appendix~\ref{app:rv-clt}) would bridge the idealized Gaussian block model with joint high-frequency realized-measure asymptotics.

\appendix
\section{High-Frequency Motivation for the Logarithmic Block Observation Model}
\label{app:rv-clt}

In this appendix, we suppress the asset index ($d=1$) and record the formal fixed-$\varepsilon$ high-frequency expansion motivating the logarithmic block observation channel \eqref{eq:realized-measure-observation}. On a filtered probability space $(\Omega,\cF,(\cF_t)_{t\ge 0},\Prob)$, consider a continuous-time asset price $S_t$ modeled as an It\^o semimartingale:
\begin{equation}
  \dd S_t = \mu_t S_t \dd t + \sigma_t^\varepsilon S_t \dd W_t^S,
  \qquad\text{or equivalently}\qquad
  \dd \log S_t = \left(\mu_t - \frac{1}{2}(\sigma_t^\varepsilon)^2\right)\dd t + \sigma_t^\varepsilon \dd W_t^S,
  \label{eq:asset-price-sde}
\end{equation}
where $\mu_t$ is a locally bounded drift process, $W_t^S$ is a standard Brownian motion, and $V_t^\varepsilon := (\sigma_t^\varepsilon)^2$ is the latent instantaneous spot variance. Sampling $S_t$ on the intra-block grid $t_{k,j} = t_{k-1} + jh$ ($j = 1,\dots,M$) with mesh $h := \delta/M$ yields the log-returns $r_{k,j} := \log S_{t_{k,j}} - \log S_{t_{k,j-1}}$ and the realized variance estimator \cite{AndersenBollerslev1998}
\begin{equation}
  \mathrm{RV}_k := \sum_{j=1}^M r_{k,j}^2,
  \qquad
  \mathrm{RV}_k \xrightarrow{\Prob} \mathrm{IV}_k := \int_{I_k} V_t^\varepsilon \dd t = \delta A_k^\varepsilon
  \quad (M \to \infty).
  \label{eq:realized-variance-def}
\end{equation}
Writing $\alpha_t := \mu_t - \frac{1}{2}V_t^\varepsilon$, a local It\^o expansion of the return increments yields
\[
  r_{k,j} = \alpha_{t_{k,j-1}} h + \sqrt{V_{t_{k,j-1}}^\varepsilon}\,\Delta W_{k,j}^S + o_{L^2}(h^{1/2}),
  \qquad \Delta W_{k,j}^S \sim \mathcal{N}(0,h).
\]
Consequently, up to terms negligible as $M \to \infty$,
\[
  \mathrm{RV}_k - \mathrm{IV}_k
  \simeq \sum_{j=1}^M V_{t_{k,j-1}}^\varepsilon\bigl((\Delta W_{k,j}^S)^2 - h\bigr),
\]
a martingale difference array whose conditional quadratic variation satisfies
\[
  2\sum_{j=1}^M \bigl(V_{t_{k,j-1}}^\varepsilon\bigr)^2 h^2
  = 2h \int_{I_k} \bigl(V_s^\varepsilon\bigr)^2 \dd s + o(h).
\]
Under standard continuous-semimartingale regularity and finite integrated quarticity \cite{BarndorffNielsenShephard2002,JacodProtter2012}, the realized-variance central limit theorem gives, for fixed $\varepsilon$ and $\delta$,
\begin{equation}
  \frac{\mathrm{RV}_k - \mathrm{IV}_k}{\sqrt{(2\delta/M)\int_{I_k}(V_s^\varepsilon)^2\dd s}}
  \Longrightarrow \mathcal{N}(0,1)
  \qquad (M \to \infty).
  \label{eq:intro-rv-clt}
\end{equation}
Since $\mathrm{IV}_k = \delta A_k^\varepsilon$ and $H_k^\varepsilon = \log A_k^\varepsilon$, the delta method yields
\[
  \log\!\left(\frac{\mathrm{RV}_k}{\delta}\right) - H_k^\varepsilon
  = \log\!\left(\frac{\mathrm{RV}_k}{\mathrm{IV}_k}\right)
  = \frac{\mathrm{RV}_k - \mathrm{IV}_k}{\mathrm{IV}_k}
    + o\!\left(\frac{|\mathrm{RV}_k - \mathrm{IV}_k|}{\mathrm{IV}_k}\right),
\]
so the logarithmic estimation error is asymptotically mixed Gaussian with path-dependent conditional variance
\begin{equation}
  R_{k,M}^{\mathrm{RV}}
  \simeq \frac{2\delta}{M}\frac{\int_{I_k}(V_s^\varepsilon)^2\dd s}{\bigl(\int_{I_k} V_s^\varepsilon \dd s\bigr)^2}.
  \label{eq:intro-rv-log-variance}
\end{equation}
Whereas $R_{k,M}^{\mathrm{RV}}$ depends on the latent path via the quarticity ratio and a joint limit in $(\varepsilon,\delta,h)$ would require a triangular-array theory resolving the fast scale ($h/\varepsilon \to 0$), the filtering analysis in Sections~\ref{sec:preliminaries}--\ref{sec:convergence} specifies a deterministic covariance $R_\delta$ and signal-independent Gaussian errors in \eqref{eq:realized-measure-observation} in order to isolate the homogenization of the nonlinear block functional $H_k^\varepsilon$.

\end{document}